\documentclass[reqno,12pt]{amsart}
\usepackage{amsmath,amsthm,amssymb,amsfonts,amscd,graphicx,xy}
\input xy
\xyoption{all}
\theoremstyle{plain} 
	\newtheorem{thm}{Theorem}[section]
	\newtheorem*{thm*}{Theorem}
	\newtheorem{cor}[thm]{Corollary}
	\newtheorem*{cor*}{Corollary}
	\newtheorem{lem}[thm]{Lemma}
	
	\newtheorem{prop}[thm]{Proposition}
	
	\newtheorem*{conj*}{Conjecture}
	
\theoremstyle{definition}
	\newtheorem{defn}[thm]{Definition}

\theoremstyle{remark}
	\newtheorem{rem}[thm]{Remark}
	
	\newtheorem*{pf}{Proof}
\numberwithin{equation}{section}

\def\CC{{\mathbb C}}

\def\PP{{\mathbb P}}
\def\QQ{{\mathbb Q}}
\def\RR{{\mathbb R}}
\def\ZZ{{\mathbb Z}}

\def\C{{\mathcal C}}
\def\D{{\mathcal D}}

\def\F{{\mathcal F}}
\def\G{{\mathcal G}}
\def\H{{\mathcal H}}
\def\I{{\mathcal I}}

\def\M{{\mathcal M}}
\def\N{{\mathcal N}}
\def\O{{\mathcal O}}

\def\T{{\mathcal T}}
\def\X{{\mathcal X}}

\def\h{{\mathfrak h}}
\def\m{{\mathfrak m}}
\def\p{\partial }

\def\ns{{\nabla}\hspace{-1.4mm}\raisebox{0.3mm}{\text{\footnotesize{\bf /}}}}

\begin{document}
\title{On the Frobenius Manifolds for Cusp Singularities}
\date{\today}
\author{Yuuki Shiraishi}
\address{Department of Mathematics, Graduate School of Science, Osaka University, 
Toyonaka Osaka, 560-0043, Japan}
\email{sm5021sy@ecs.cmc.osaka-u.ac.jp}
\author{Atsushi Takahashi}
\address{Department of Mathematics, Graduate School of Science, Osaka University, 
Toyonaka Osaka, 560-0043, Japan}
\email{takahashi@math.sci.osaka-u.ac.jp}
\begin{abstract}
We show that 
the Frobenius manifold associated to the pair of a cusp singularity and 
it's canonical primitive form is isomorphic to the one constructed from 
the Gromov--Witten theory for an orbifold projective line with at most three 
orbifold points. We also calculate the intersection form of the Frobenius manifold.
\end{abstract}
\maketitle
\section{Introduction}

The notion of primitive forms is introduce by K. Saito \cite{sa:1, S1202-Saito} in his study of period mapping associated to a deformation theory of an isolated hypersurface singularity. Roughly speaking,  a primitive form is a generalization of a differential of the first kind on an elliptic curve and the tools in order to define a primitive form are natural generalization of a polarized Hodge structure. The existence of primitive forms is proved for isolated hypersurface singularities \cite{m.saito} (for comprehensive guides, see also \cite{her:1, sabb:1}).

In the first part of the present paper, we shall investigate the canonical primitive form for
a cusp singularity: 
Let $A$ be a triplet $(a_1,a_2,a_3)$ of positive integers such that $a_1\le a_2\le a_3$.
Set $\mu_A=a_1+a_2+a_3-1$ and $\chi_A:=1/a_1+1/a_2+1/a_3-1$.
We can associate to $A$ the polynomial $f_A({\bf x})\in\CC[x_1,x_2,x_3]$ given as
\begin{equation}
f_A({\bf x}):=x_1^{a_1}+x_2^{a_2}+x_3^{a_3}-q^{-1}\cdot x_1 x_2 x_3
\end{equation}
for some $q\in\CC\backslash\{0\}$, which we shall call a cusp polynomial of type $A$. In particular, it defines a cusp singularity if $\chi_A<0$.
It was claimed by K. Saito \cite{S1202-Saito} and is proven by M.~Saito \cite{m.saito} that, for the universal unfolding of 
a cusp singularity, there exists a unique primitive form up to a constant factor with the minimal exponent $r=1$, whose associated exponents are given by the mixed Hodge structure of the cusp singularity. 
First, we shall give its local expression in order to show the mirror isomorphism:
\begin{thm*}[Theorem \ref{primitive form}]
Assume that $\chi_A<0$.
There exists a unique primitive form $\zeta_A$ for the tuple $(\H_{F_A}^{(0)},\nabla, K_{F_A})$ 
with the minimal exponent $r=1$ such that 
\begin{equation}
\zeta_A|_{{\bf s}={\bf 0}}=[s_{\mu_A}^{-1}dx_1 \wedge dx_2 \wedge dx_3].
\end{equation}
\end{thm*}
Once we obtain the primitive form, we can construct the Frobenius structure on
the deformation space of the isolated hypersurface singularity 
by well--known construction \cite{st:1}. 
Next, we shall show that the Frobenius structure associated to 
the pair $(f_A, \zeta_A)$ given above satisfies the conditions in Theorem 3.1 of \cite{ist:1}:
\begin{thm*}[Theorem~\ref{thm:satisfies IST1}]
Assume that $\chi_A<0$.
For the Frobenius structure of rank $\mu_A$ and dimension one constructed from the pair $(f_A,\zeta_A)$, 
there exist flat coordinates $t_1$, $t_{1,1},\dots ,t_{3,a_3-1}$, $t_{\mu_A}$ satisfying the following conditions$:$
\begin{enumerate}
\item 
The unit vector field $e$ and the Euler vector field $E$ are given by
\[
e=\frac{\p}{\p t_1},\ E=t_1\frac{\p}{\p t_1}+\sum_{i=1}^3\sum_{j=1}^{a_i-1}\frac{a_i-j}{a_i}t_{i,j}\frac{\p}{\p t_{i,j}}
+\chi_A\frac{\p}{\p t_{\mu_A}}.
\]
\item 
The non-degenerate symmetric bilinear form $\eta$ on $\T_M$ defined by 
\begin{equation*}
\eta(\delta,\delta'):=K_{F_A}^{(0)}(u\nabla_\delta\zeta_A,u\nabla_{\delta'}\zeta),\quad \delta,\delta'\in\T_M
\end{equation*}
satisfies
\begin{align*}
&\ \eta\left(\frac{\p}{\p t_1}, \frac{\p}{\p t_{\mu_A}}\right)=
\eta\left(\frac{\p}{\p t_{\mu_A}}, \frac{\p}{\p t_1}\right)=1,\\ 
&\ \eta\left(\frac{\p}{\p t_{i_1,j_1}}, \frac{\p}{\p t_{i_2,j_2}}\right)=
\begin{cases}
\frac{1}{a_{i_1}}\quad i_1=i_2\text{ and }j_2=a_{i_1}-j_1,\\
0 \quad \text{otherwise}.
\end{cases}
\end{align*}
\item 
The Frobenius potential $\F_{f_A,\zeta_A}$ satisfies $E\F_{f_A,\zeta_A}|_{t_{1}=0}=2\F_{f_A,\zeta_A}|_{t_{1}=0}$,
\[
\left.\F_{f_A,\zeta_A}\right|_{t_1=0}\in\CC\left[[t_{1,1}, \dots, t_{1,a_1-1},
t_{2,1}, \dots, t_{2,a_2-1},t_{3,1}, \dots, t_{3,a_3-1},e^{t_{\mu_A}}]\right].
\]
\item The restriction of the Frobenius potential $\F_{f_A,\zeta_A}$ to 
the submanifold $\{t_1=e^{t_{\mu_A}}=0\}$ is given as 
\[
\F_{f_A,\zeta_A}|_{t_1=e^{t_{\mu_A}}=0}=\G^{(1)}+\G^{(2)}+\G^{(3)}, 
\]
where $\G^{(i)}\in \CC[t_{i,1},\dots, t_{i,a_i-1}]$, $i=1,2,3$.
\item 
In the frame $\frac{\p}{\p t_1}$, $\frac{\p}{\p t_{1,1}},\dots, 
\frac{\p}{\p t_{3,a_3-1}}$, $\frac{\p}{\p t_{\mu_A}}$ of $\T_M$,
the product $\circ$ can be extended to the limit $t_1=t_{1,1}=\dots=t_{3,a_3-1}=e^{t_{\mu_A}}=0$.
The $\CC$-algebra obtained in this limit is isomorphic to 
\[
\CC[x_1,x_2,x_3]\left/\left(x_1x_2,x_2x_3,x_3x_1,a_1x_1^{a_1}-a_2x_2^{a_2},
a_2x_2^{a_2}-a_3x_3^{a_3}\right)\right.,
\]
where $\p/\p t_{1,1},\p/\p t_{2,1},\p/\p t_{3,1}$ are mapped to
$x_1,x_2,x_3$, respectively.
\item The term $t_{1,1}t_{2,1}t_{3,1}e^{t_{\mu_A}}$
occurs with the coefficient $1$ in $\F_{f_A,\zeta_A}$. 
\end{enumerate}
\end{thm*}
As a consequence, 
we obtain the mirror isomorphism as Frobenius manifolds:
\begin{thm*}[Theorem~\ref{cor:classical mirror}]
Assume that $\chi_A<0$.
There exists an isomorphism of Frobenius manifolds between the one constructed from the Gromov--Witten theory 
for $\PP^{1}_A$ and the one constructed from the pair $(f_A,\zeta_A)$. 
\end{thm*}
For the cases that $\chi_A>0$ by \cite{mt:1, r:1, ist:2} and the cases that $\chi_A=0$ by \cite{sat:1},
same statements as Theorem ~\ref{cor:classical mirror} are shown.
Therefore, combining them with Theorem ~\ref{cor:classical mirror}, it is shown that,
for arbitary triplet of positive integers $A$, there exists the classical mirror symmetry between
the orbifold projective line with at most three orbifold points $\PP^{1}_A$ and the pair of the cusp polynomial $f_A$ and the primitive form $\zeta_A$ associated to it.
In the last part of the present paper, we shall investigate the period mapping of the primitive form $\zeta_A$ for the case that $\chi_A\ne 0$.
Here $\zeta_A$ is the one given in Theorem \ref{primitive form} for $\chi_A<0$ 
and the one in Theorem 3.1 in \cite{ist:2} for $\chi_A>0$. Let $T_A$ be the following
Coxeter--Dynkin diagram:
\begin{equation}
\begin{xy}
\xygraph{
    \bullet ([]!{+(0,-.3)} {(1,a_1-1)}) - [r] \cdots - [r]
    \bullet ([]!{+(0,-.3)} {(1,1)}) - [r] 
    \bullet ([]!{+(.3,-.3)} {1}) (
        - [d] \bullet ([]!{+(.5,0)} {(2,1)})( 
        - [d] \vdots - [d] \bullet ([]!{+(.8,0)} {(2,a_2-1)})),
        - [r] \bullet ([]!{+(0,-.3)} {(3,1)})
        - [r] \cdots - [r] \bullet ([]!{+(0,-.3)} {(3,a_3-1)})
)}
\end{xy}
\end{equation}
and $\h_A$ the complexified Cartan subalgebra of the Kac--Moody Lie algebra associated to $T_A$
(in particular, which is the simple Lie algebra for the case that $\chi_A>0$).
Denote by $\alpha_1,\dots, \alpha_{(3, a_3-1)}\in \h_A^{*}:={\rm Hom}_{\CC}(\h_A,\CC)$ simple roots corresponding the vertices in $T_{A}$, 
by $\alpha^{\vee}_1,\dots, \alpha^{\vee}_{(3, a_3-1)}\in \h_A$ simple coroots
and by $\left<, \right>:\h^{*}_A\otimes_{\CC} \h_A\longrightarrow \CC$ the natural pairing.

By using some isomorphisms (see Lemma \ref{Fourier-Laplace of primitve form}, Lemma \ref{isomorphism of lattice} and Lemma \ref{subspace of ah}), 
we can identify $\alpha_{\bf i}, \ {\bf i}=1, (1,1), \dots, (i,j), \dots, (3, a_{3}-1)$ with homology classes represented by vanishing cycles in each fiber $\X_{(0, {\bf s}, s_{\mu_A})}\subset \CC^{3}$ (see Subsection \ref{preliminary period}) over a general point $(0, {\bf s}, s_{\mu_A})$.
Under this situation, we shall consider periods of the holomorphic $2$-form $\tilde{\zeta_A}$ on $\X_{(0, {\bf s}, s_{\mu_A})}$ induced from the primitive form $\zeta_A$ 
(see Subsection 5.2, \eqref{zetatilde}).
Then we can calculate the intersection form $I_{(f_A,\zeta_A)}$ of $M:=M_{(f_A,\zeta_A)}$:

\begin{thm*}[Theorem \ref{calc: intersection form 2}]
Consider the periods
\begin{equation}
x_{\bf i}:=\frac{1}{(2\pi\sqrt{-1})^{2}}\int_{\gamma_{\bf i}({\bf s},s_{\mu_A})} \tilde{\zeta_A}, \quad {\bf i}=1, (1,1), \dots, (i,j), \dots, (3, a_{3}-1),
\end{equation}
where $\gamma_{\bf i}({\bf s},s_{\mu_A})$ is the horizontal family of homology classes 
in $\displaystyle\bigcup_{\substack{(0, {\bf s}, s_{\mu_A}) \in \check{M} \backslash \check{\D}}} H_{2}(\X_{(0, {\bf s}, s_{\mu_A})},\CC)$ identified with $\alpha_{\bf i}$, and the function
\begin{equation}
x_{\mu_A}:=\frac{1}{2\pi \sqrt{-1}} t_{\mu_A}=\frac{1}{2\pi \sqrt{-1}} \log s_{\mu_A}.
\end{equation}
They define flat coordinates with respect to $I_{(f_A,\zeta_A)}$ on the monodromy covering space of $M\backslash \D$.
Moreover, one has
\begin{subequations}
\begin{align}
\lefteqn{I_{{(f_A,\zeta_A)}} (dx_{\bf i}, dx_{\bf j})
=\frac{-1}{(2\pi \sqrt{-1})^2}\left<\alpha_{\bf i},\alpha_{\bf j}^{\vee}\right>,}\\
&I_{{(f_A,\zeta_A)}} (dx_{\mu_A}, dx_{\bf i})
=I_{{(f_A,\zeta_A)}} (dx_{\bf i}, dx_{\mu_A})=0,\\
&I_{{(f_A,\zeta_A)}}(dx_{\mu_A},dx_{\mu_A})
=\displaystyle \frac{1}{(2\pi \sqrt{-1})^{2}}\chi_A.
\end{align}
\end{subequations}
\end{thm*}

On the other hand, Dubrovin--Zhang \cite{dz:1} constructed the Frobenius structure on a quotient space by an extended affine
Weyl group $\widehat{W}_A$ for $\widehat{\h}_A:=\h_A\times \CC$: 
\begin{thm*}[\cite{dz:1}]
Assume that $\chi_A>0$.
There exists a unique {\rm Frobenius} structure of rank $\mu_A$ and dimension one on 
$M_{\widehat{W}_A}:=\widehat{\h}_A/\widehat{W}_A$ with flat coordinates 
$t_{1},t_{1,1},\dots, t_{i,j}, \dots, t_{3,a_{3}-1}, t_{\mu_A}:=(2\pi\sqrt{-1})x_{\mu_A}$ such that
\begin{equation}
e=\frac{\p}{\p t_1}, \ \ E=t_1\frac{\p}{\p t_1}+\sum_{i=1}^3\sum_{j=1}^{a_i-1}\frac{a_i-j}
{a_i}t_{i,j}\frac{\p}{\p t_{i,j}}
+\chi_A\frac{\p}{\p t_{\mu_A}},
\end{equation}
and the intersection form $I_{\widehat{W}_A}$ is given by
\begin{subequations}
\begin{align}
\lefteqn{
I_{\widehat{W}_A}(\alpha_{\bf i}, \alpha_{\bf j})=\frac{-1}{(2\pi\sqrt{-1})^{2}}
\left<\alpha_{\bf i}, \alpha^{\vee}_{\bf j}\right>, \quad {\bf i,j}=1, (1,1),\dots,(3,a_{3}-1),}\\
&I_{\widehat{W}_A}(\alpha_{\bf i},dx_{\mu_A})=I_{\widehat{W}_A}(dx_{\mu_A},\alpha_{\bf i})
=0,\quad {\bf i}=1, (1,1),\dots,(3,a_{3}-1),\\
&I_{\widehat{W}_A}(dx_{\mu_A},dx_{\mu_A})
=\frac{1}{(2\pi\sqrt{-1})^{2}}\chi_A,
\end{align}
\end{subequations}
where we identify the cotangent space of $M_{\widehat{W}_A}$
with $\h^{*}_A\oplus\CC dx_{\mu_A}$.
\end{thm*}

The isomorphism as Frobenius manifolds between $M_{(f_A, \zeta_A)}$ and $M_{\widehat{W}_A}$
is obtained in \cite{dz:1} and \cite{mt:1} for the case that $A=(1, a_2, a_3)$ 
and \cite{r:1} for left cases such that $\chi_A>0$. However, 
systematical proof is not found in any literature. By contrast we can give a natural proof
by using Theorem \ref{calc: intersection form 2} and the following theorem which might be known to experts:
\begin{thm*}[Theorem \ref{second}]
A Frobenius manifold $M$ of rank $\mu_A$ and dimension one with the 
following $e$ and $E$
is uniquely determined by the intersection form $I_M$ $:$
\begin{equation}
e=\frac{\p}{\p t_1},\ E=t_1\frac{\p}{\p t_1}+\sum_{i=1}^3\sum_{j=1}^{a_i-1}\frac{a_i-j}
{a_i}t_{i,j}\frac{\p}{\p t_{i,j}}
+\chi_A\frac{\p}{\p t_{\mu_A}}.
\end{equation}
\end{thm*}

As a consequence, for the case that $\chi_A>0$,
we can obtain an isomorphism as Frobenius manifolds between $M_{(f_A,\zeta_A)}$ and $M_{\widehat{W}_A}$
in a natural way:
\begin{cor*}[Corollary \ref{cor:sing and weyl}]
Assume that $\chi_A>0$.
There exists an isomorphism of Frobenius manifolds between 
the one constructed from the invariant thory of extended affine Weyl group 
$\widehat{W}_A$ and the one constructed from the pair $(f_A,\zeta_A)$. 
\end{cor*}
\bigskip
\noindent
{\it Acknowledgement}\\
\indent
The first named author is deeply grateful to Claus Hertling and 
Christian Sevenheck for their valuable discussions and encouragement.  
He is supported by the JSPS International Training Program (ITP).
The second named author is supported by JSPS KAKENHI Grant Number
24684005.
\section{Notations and Terminologies}\label{sec:N and T}
Let $A$ be a triplet $(a_1,a_2,a_3)$ of positive integers such that $a_1\le a_2\le a_3$.
Set 
\begin{equation}
\mu_A:=a_1+a_2+a_3-1
\end{equation}
and
\begin{equation}
\chi_A:=\frac{1}{a_1}+\frac{1}{a_2}+\frac{1}{a_3}-1.
\end{equation}
\subsection{Universal unfolding of cusp singularities}
\begin{defn}
A polynomial $f_A({\bf x})\in\CC[x_1,x_2,x_3]$ given as
\begin{equation}
f_A({\bf x}):=x_1^{a_1}+x_2^{a_2}+x_3^{a_3}-q^{-1}\cdot x_1 x_2 x_3
\end{equation}
for some $q\in\CC\backslash\{0\}$ is called 
the {\it cusp polynomial of type $A$}. 
In particular, it defines a {\it cusp singurality} if $\chi_A<0$. 
\end{defn}
Let $f_A=f_A({\bf x})$ be a cusp singularity of type $A$ such that $\chi_A<0$.
Then we regard $f_A$ as a {\it holomorphic function} defined on 
a neighborhood of the origin in $\CC^{3}$. 
\begin{prop}
The germ of the hypersurface $X_0:=\{f_A({\bf x})=0\}\subset (\CC^{3},0)$ 
at the origin has at most an isolated singular point at the origin. 
In particular, we have  
\begin{equation}
\O_{\CC^{3},0}\left/\Big(\frac{\p f_A}{\p x_1}, \frac{\p f_A}{\p x_2}, \frac{\p f_A}{\p x_3}\Big)\right.
\end{equation}
is a $\CC$-vector space of rank $\mu_A$.
\end{prop}
\begin{pf}
This is a well-known fact. Note that we can choose a basis of the vector space as
\begin{equation}
1,\ x_i^j\ (i=1,2,3,j=1,\dots, a_i-1),\ x_1x_2x_3.
\end{equation} 
\qed
\end{pf}
We can consider the {\it universal unfolding} of $f_A$.
\begin{defn}
A holomorphic function $F_A({\bf x};{\bf s}, s_{\mu_A})$ defined on a neighborhood of $({\bf 0};{\bf 0},q)$ of 
$\CC^{3}\!\times\! (\CC^{\mu_A-1}\!\times\! \CC\backslash\{0\})$ is given as follows$;$ 
\begin{equation}
F_A({\bf x};{\bf s},s_{\mu_A}):=x_1^{a_1}+x_2^{a_2}+x_3^{a_3}-s_{\mu_A}^{-1}\cdot x_1 x_2 x_3
+s_1\cdot 1+\sum_{i=1}^{3} \sum_{j=1}^{a_i-1}s_{i,j}\cdot x_i^j.
\end{equation}
\end{defn}
In order to get a flat family of hypersurfaces where the fiber over $({\bf 0};{\bf 0},q)$ is isomorphic 
to the singularity $X_0$ and generic fibers are smoothings of $X_0$ containing vanishing cycles to the singularity,
we have to avoid ``cycles from infinity".
Therefore, we shall choose suitably a domain and a range of $F_A$ and an open set 
$\CC^{\mu_A-1}\!\times\!\CC\backslash\{0\}$ in the following way.
Denote by
\begin{equation}
p:\CC^{3}\!\times\!(\CC^{\mu_A-1}\!\times\!\CC)\longrightarrow \CC^{\mu_A-1}\!\times\! \CC,
\quad ({\bf x};{\bf s},s_{\mu_A})\mapsto ({\bf s},s_{\mu_A})
\end{equation}
the projection map.
Fix Euclidean norms $\| \cdot \|$ 
on $\CC^{3}\!\times\!(\CC^{\mu_A-1}\!\times\!\CC)$, $\CC$ 
and $\CC^{\mu_A-1}\!\times\!\CC$.
For positive real numbers $r,\delta$ and $\epsilon$, put
\begin{subequations}
\begin{equation}
\CC_{\delta }:=\{w\in \CC~|~ \| w\| <\delta\},
\end{equation}
\begin{equation}
M_{\epsilon }:=\{({\bf s},s_{\mu_A})\in \CC^{\mu_A-1}\!\times\!\CC\backslash\{0\}~|~
\|({\bf s},s_{\mu_A})\|<\epsilon\},
\end{equation}
\begin{equation}
\overline{M}_{\epsilon }:=\{({\bf s},s_{\mu_A})\in \CC^{\mu_A-1}\!\times\!\CC~|~
\|({\bf s},s_{\mu_A})\|<\epsilon\},
\end{equation}
\begin{multline}
\X_{r,\delta ,\epsilon }
:=\{({\bf x};{\bf s},s_{\mu_A})\in \CC^{3}\!\times\!(\CC^{\mu_A-1}\!\times\!\CC\backslash\{0\})~|~\\
\|({\bf x};{\bf s},s_{\mu_A})\| <r, \| F_A({\bf x};{\bf s},s_{\mu_A})\| <\delta\} 
\cap p^{-1}(M_{\epsilon}),
\end{multline}
\begin{multline}
\overline{\X}_{r,\delta ,\epsilon }
:=\{({\bf x};{\bf s},s_{\mu_A})\in \CC^{3}\!\times\!(\CC^{\mu_A-1}\!\times\!\CC)~|~\\
\|({\bf x};{\bf s},s_{\mu_A})\| <r, \| F_A({\bf x};{\bf s},s_{\mu_A})\| <\delta\} 
\cap p^{-1}(M_{\epsilon}).
\end{multline}
\end{subequations}
For the choice $1\gg r\gg\delta \gg\epsilon> 0$ of the radius, we have 
\begin{enumerate}
\item $p:\X_{r,\delta ,\epsilon }\to M_{\epsilon }$ is a smooth Stein map, 
which is topologically acyclic.
\item The holomorphic map 
\begin{equation}
\varphi:\X_{r,\delta ,\epsilon}\longrightarrow \CC_{\delta}\times M_{\epsilon},
\quad ({\bf x};{\bf s}, s_{\mu_A})\mapsto (w;{\bf s},s_{\mu_A}):=
(F_A({\bf x};{\bf s},s_{\mu_A});{\bf s},s_{\mu_A})
\end{equation} 
is a flat Stein map, whose fibers are smooth and transverse to $\partial X_{r,\delta.\epsilon}$ 
at each point of $\partial \X_{r,\delta.\epsilon}$. 
The fiber $\varphi^{-1}({\bf 0};{\bf 0},q)\cap \X_{r,\delta.\epsilon}$ over $({\bf 0};{\bf 0},q)$ has 
only an isolated singular point isomorphic to $X_0$ and is contractible to a point.
\end{enumerate}
We fix the constants $r,\delta$ and $\epsilon$ once for all and 
put $\X:=\X_{r,\delta,\epsilon}$, $\overline{\X}:=\overline{\X}_{r,\delta,\epsilon}$, 
$M:=M_{\epsilon}$ and $\overline{M}:=\overline{M}_\epsilon$.
We also introduce the following notations for the later convenience:
\begin{defn}
Put
\begin{subequations}
\begin{equation}
\X^s:=\{ ({\bf x};{\bf s},s_{\mu_A})\in\X~|~{\bf s}={\bf 0}\},
\end{equation}
\begin{equation}
M^s:=\{ ({\bf s},s_{\mu_A})\in M~|~{\bf s}={\bf 0}\},
\end{equation}
\begin{equation}
F_A^s:=F_A({\bf x};{\bf 0},s_{\mu_A}).
\end{equation}
\end{subequations}
\end{defn}
\begin{rem}
$M^s$ will correspond by mirror symmetry to the ``small quantum cohomology subspace", 
that is the reason of the letter ``s".
\end{rem}
Set 
\begin{equation}
\O_\C:=\O_{\X}
\left/\Big(\frac{\p F_A}{\p x_1},\frac{\p F_A}{\p x_2},\frac{\p F_A}{\p x_3}\Big)\right..
\end{equation}
\begin{prop}\label{unfolding}
The function $F_A({\bf x};{\bf s},s_{\mu_A})$ satisfies the following conditions$:$
\begin{enumerate}
\item $F_A({\bf x};{\bf 0},q)=f_A({\bf x}).$
\item The $\O_{M}$-homomorphism 
$\rho$ called the Kodaira--Spencer map defined as 
\begin{equation}\label{universal}
\rho: \T_{M}\longrightarrow p_*\O_\C,\quad \delta\mapsto \widehat{\delta} F_A,
\end{equation}
is an isomorphism, where $\widehat{\delta}$ is a lifting on $\X$ 
of a vector field $\delta$ with respect to the projection $p$.
\end{enumerate}
\end{prop}
\begin{pf}
This is also a well-known fact.
\qed
\end{pf}
Note that the tangent bundle $\T_M$ naturally obtains an $\O_M$-algebra structure. 
\begin{defn}
We shall denote by $\circ$ the induced product structure on $\T_M$ by the $\O_M$-isomorphism~\eqref{universal}. 
Namely, for $\delta,\delta'\in \T_M$, we have
\begin{equation}
\widehat{(\delta\circ\delta')}F_A = 
\widehat{\delta}F_A\cdot\widehat{\delta'}F_A\ \text{in}\ p_*\O_\C.
\end{equation}
\end{defn}
\begin{defn}
Let $\iota$ be the natural embedding
\begin{equation}
\iota:M\hookrightarrow \overline{M}, \quad
({\bf s}, s_{\mu_A}) \mapsto ({\bf s}, s_{\mu_A}).
\end{equation}
Define by $M_\infty$ the divisor 
\begin{equation}
M_\infty :=\overline{M}\backslash\iota(M):=\{ ({\bf s},s_{\mu_A})\in\overline{M}~|~s_{\mu_A}=0\}.
\end{equation}
\end{defn}
Set 
\begin{equation}\label{eq:Jacobi-ring extension}
\O_{\overline{\C}}:=\O_{\overline{\X}}
\left/\Big(s_{\mu_A}\frac{\p F_A}{\p x_1},s_{\mu_A}\frac{\p F_A}{\p x_2},s_{\mu_A}\frac{\p F_A}{\p x_3},
H_{1}({\bf x}, {\bf s}), H_{2}({\bf x}, {\bf s})
\Big)\right..
\end{equation}
where
\begin{equation}
H_{i}({\bf x}, {\bf s}):=a_{i}x_{i}^{a_i}-a_{i+1}x_{i+1}^{a_{i+1}}+
\sum_{j=1}^{a_i-1}j\cdot s_{i,j}\cdot x_i^j-\sum_{j=1}^{a_{i+1}-1}j\cdot s_{i+1,j}\cdot x_{i+1}^j,\ i=1,2.
\end{equation}
\begin{prop}
As an $\O_{\overline{M}}$-module, ${p}_*\O_{\overline{\C}}$ is free of rank $\mu_A$.
As an $\O_{\overline{M}}$-algebra, it is isomorphic to the $\O_{\overline{M}}$-subalgebra of $\iota_*\left(p_*\O_\C\right)$ 
whose $\O_{\overline{M}}$-basis is given by the set of residue classes of monomials 
\begin{equation}
\left\{1,\ x_i^j\  (i=1,2,3,\ j=1,\dots, a_i-1),\ s_{\mu_A}^{-1}x_1x_2x_3\right\}.
\end{equation}
In particular, we have the isomorphism of $\O_M$-modules
\begin{equation}
\iota^*\left({p}_*\O_{\overline{\C}}\right)\simeq p_*\O_\C.
\end{equation}
Namely, ${p}_*\O_{\overline{\C}}$ is an extension of $p_*\O_\C$ as a free $\O_{\overline{M}}$-module.
\end{prop}
\begin{pf}
It is almost obvious that $\iota^*\left({p}_*\O_{\overline{\C}}\right)\simeq p_*\O_\C$.
\begin{lem}\label{lem:limit of Jac}
We have
\begin{equation}
{p}_*\O_{\overline{\C}}\left/s_{\mu_A}{p}_*\O_{\overline{\C}}\right.
\simeq\O_{M_\infty}[x_1,x_2,x_3]\left/\Big(x_2x_3,x_3x_1,x_1x_2,H_{1}({\bf x}, {\bf s}), H_{2}({\bf x}, {\bf s})\Big)\right.,
\end{equation}
is a free $\O_{M_\infty}$-module of rank $\mu_A$. In particular, we have
\begin{equation}
{p}_*\O_{\overline{\C}}\left/\m_{({\bf 0},0)}{p}_*\O_{\overline{\C}}\right.
\simeq\CC[x_1,x_2,x_3]\left/\Big(x_2x_3,x_3x_1,x_1x_2,a_1x_1^{a_1}-a_2x_2^{a_2},a_2x_2^{a_2}-a_3x_3^{a_3}\Big)\right.,
\end{equation}
where $\m_{({\bf 0},0)}$ is the maximal ideal of $\O_{\overline{M}}$ corresponding to the point $({\bf 0},0)$.
\end{lem}
\begin{pf}
Some elementary calculations yield the statement.
\qed
\end{pf}
By this lemma, we see that ${p}_*\O_{\overline{\C}}$ is free of rank $\mu_A$ as an $\O_{\overline{M}}$-module.
The rest follows from the following equality in $\iota_*\left(p_*\O_\C\right):$
\begin{equation}
\left[s_{\mu_A}^{-1}x_1x_2x_3\right]=\left[a_ix_i^{a_{i}}+\sum_{j=1}^{a_i-1}j\cdot s_{i,j}\cdot x_i^j\right],\quad i=1,2,3.
\end{equation}
\qed
\end{pf}
Therefore, the analytic subset $\overline{\C}$ in $\overline{\X}$ defined by \eqref{eq:Jacobi-ring extension}
is a closure of the relative critical set $\C$, 
which is flat and finite over $\overline{M}$.
\begin{prop}\label{prop:limit of Jac}
Denote by $\T_{\overline{M}}\left(-\log M_\infty\right)$ the sheaf of holomorphic vector fields on $\overline{M}$ 
with logarithmic zeros along the divisor $M_\infty$ defined as 
\begin{equation}
\T_{\overline{M}}\left(-\log M_\infty\right):=\{\delta\in\T_{\overline{M}}~|~\delta s_{\mu_A}\in s_{\mu_A}\cdot \O_{\overline{M}}\}.
\end{equation}
Then, $\T_{\overline{M}}\left(-\log M_\infty\right)$ is a free $\O_{\overline{M}}$-module of rank $\mu_A$ 
and we have the isomorphism of $\O_M$-modules
\begin{equation}
\iota^*\left(\T_{\overline{M}}\left(-\log M_\infty\right)\right)\simeq \T_M.
\end{equation}
Namely, $\T_{\overline{M}}\left(-\log M_\infty\right)$ is an extension of $\T_M$ as a free $\O_{\overline{M}}$-module.
\end{prop}
\begin{pf}
The statement follows from the fact that 
\begin{equation}
\T_{\overline{M}}\left(-\log M_\infty\right)\simeq \O_{\overline{M}}\frac{\p}{\p s_1}
\bigoplus_{\substack{
1\le i\le 3, \\
1\le j\le a_{i}-1}} 
\O_{\overline{M}}\frac{\p}{\p s_{i,j}}
\bigoplus \O_{\overline{M}}s_{\mu_A}\frac{\p}{\p s_{\mu_A}}.
\end{equation}
\qed
\end{pf}
\begin{prop}\label{prop:KS-extension} 
The Kodaira--Spencer map $\rho$ induces the $\O_{\overline{M}}$-isomorphism $\overline{\rho}$
\begin{equation}
\overline{\rho}:\T_{\overline{M}}\left(-\log M_\infty\right)\longrightarrow {p}_*\O_{\overline{\C}}.
\end{equation}
\end{prop}
\begin{pf}
The statement easily follows since the Kodaira--Spencer map $\rho$ induces the $\O_{\overline{M}}$-isomorphism 
between $\T_{\overline{M}}\left(-\log M_\infty\right)$ and the free 
$\O_{\overline{M}}$-submodule of $\iota_*\left(p_*\O_\C\right)$ spanned by the residue classes of monomials 
\begin{equation}
1,\ x_i^j\  (i=1,2,3,\ j=1,\dots, a_i-1),\ s_{\mu_A}^{-1}x_1x_2x_3.
\end{equation}
\qed
\end{pf}
\subsection{Primitive vector field and Euler vector field}
\begin{defn}
The vector field $e$ and $E$ on $M$ corresponding to the unit 
$1$ and $F$ by the $\O_M$-isomorphism~\eqref{universal} 
is called the {\it primitive vector field} and the {\it Euler vector field}, respectively. 
That is, 
\begin{equation}
\widehat{e}F_A=1\ \text{and}\ \widehat{E}F_A=F_A\ 
\text{in}\ p_*\O_\C.
\end{equation}
\end{defn}
\begin{prop}
The primitive vector field $e$ and the Euler vector field $E$ on $M$ are given by 
\begin{equation}
e=\frac{\p}{\p s_1},\quad E=s_1\frac{\p }{\p s_1}+\sum_{i=1}^3\sum_{j=1}^{a_i-1}\frac{a_i-j}{a_i}s_{i,j}\frac{\p}{s_{i,j}}+
\chi_A s_{\mu_A}\frac{\p}{\p s_{\mu_A}}.
\end{equation}
\end{prop}
\begin{pf}
Easy calculation yields the statement.
\qed
\end{pf}
\begin{lem}
We have the ``Euler's identity"$:$
\begin{equation}
F_A=EF_A+\sum_{i=1}^3\frac{1}{a_i}x_i\frac{\p F_A}{\p x_i}.
\end{equation}
\end{lem}
\begin{pf}
Easy calculation yields the statement.
\qed
\end{pf}
\begin{defn}
An element $g\in\O_M$ is of {\it degree} $l$ for some $k\in \QQ$ if 
it satisfies the equation $Eg=lg$ where $E$ is the Euler vector field. 
If $Eg=lg$, then $l$ is denoted by $\deg(g)$.
\end{defn}
\subsection{Filtered de Rham cohomology}
For any non-negative integer $i$, denote by $\Omega^i_{\X/M}$ the sheaf of relative holomorphic differential 
$i$-forms with respect to the projection $p:\X\longrightarrow M$.
\begin{defn}
Define $\Omega_{F_A}$ as
\begin{equation}
\Omega_{F_A}:=p_{*}\Omega^{3}_{\X/M}/dF_{A}\wedge p_{*}\Omega^{2}_{\X/M}.
\end{equation}
\end{defn}
\begin{prop}\label{omega F}
$\Omega_{F_A}$ is a free $p_*\O_\C$-module of rank one and hence a free $\O_M$-module of rank $\mu_A$.
\end{prop}
\begin{pf}
A $p_*\O_\C$-free base of $\Omega_{F_A}$ can be chosen as $[\omega]$, $\omega:=dx_1\wedge dx_2\wedge dx_3$,
and hence an $\O_M$-free basis of $\Omega_{F_A}$ can be chosen as 
\begin{equation}
\left\{[\omega],\ [x_i^j\omega]\  (i=1,2,3,j=1,\dots, a_i-1),\ [x_1x_2x_3\omega]\right\}.
\end{equation}
\qed
\end{pf}
\begin{defn}
We set 
\begin{equation}
\H_{F_A}:=\RR^3 p_{*}(\Omega^{\bullet}_{\X/M}\otimes_{\O_{M}}\O_M((u)),
u d_{\X/M}+dF_A\wedge)
\end{equation}
and call it the {\it filtered de Rham cohomology group} of the universal unfolding $F_A$.
\end{defn}
For any $k\in\ZZ$, put
\begin{equation}
\H_{F_A}^{(-k)}:=\RR^3 p_{*}(\Omega^{\bullet}_{\X/M}\otimes_{\O_{M}}\O_M[[u]]u^k,
u d_{\X/M}+dF_A\wedge).
\end{equation}
Obviously, one has an $\O_M[[u]]$-isomorphism for all $i\in\ZZ$
\begin{equation}
\H_{F_A}^{(0)}\simeq \H_{F_A}^{(-k)},\quad \omega\mapsto u^k\omega,
\end{equation}
and an $\O_M((u))$-isomorphism 
\begin{equation}
\H_{F_A}^{(-k)}\otimes_{\O_M[[u]]}\O_M((u))\simeq \H_{F_A}.
\end{equation}
Furthermore, $\H_{F_A}^{(-k)}$ is naturally a submodule of $\H_{F_A}$ 
so that $\{\H_{F_A}^{(-k)}\}_{i\in\ZZ}$ form an increasing and exhaustive filtration of $\H_{F_A}$:
\begin{equation}\label{filt}
\cdots \subset \H^{(-k-1)}_{F_A} \subset \H^{(-k)}_{F_A} \subset \cdots \subset 
\H^{(-1)}_{F_A}\subset \H^{(0)}_{F_A} \subset \cdots \subset \H_{F_A},  
\end{equation}
such that $\H_{F_A}$ is complete with respect to the filtration in the following sense:
\begin{equation}
\bigcup_{k\in \ZZ} \H^{(-k)}_{F_A}=\H_{F_A} \ \text{and} \ \bigcap_{k\in \ZZ} \H^{(-k)}_{F_A}=\{0\}.  
\end{equation}
\begin{prop}\label{b-lat-k}
For any $k\in\ZZ$, $\H_{F_A}^{(-k)}$ is an $\O_M[[u]]$-free module of rank $\mu_A$.
In particular, we have the following short exact sequence of $\O_M$-modules
\begin{equation}
0 \rightarrow \H^{(-1)}_{F_A}\hookrightarrow \H^{(0)}_{F_A} \xrightarrow{r^{(0)}} \Omega_{F_A}\rightarrow 0.
\end{equation}
\end{prop}
\begin{pf}
One can choose an $\O_M[[u]]$-free basis of $\H^{(0)}_{F_A}$ as 
\begin{equation}
\left\{[\omega],\ [x_i^j\omega]\  (i=1,2,3,j=1,\dots, a_i-1),\ [x_1x_2x_3\omega]\right\}
\end{equation}
where $\omega=dx_1\wedge dx_2\wedge dx_3$. The rest is clear.
\qed
\end{pf}
\begin{defn}
Define an element $\zeta_A'$ of $\Gamma(M,\H_{F_A}^{(0)})$ as 
\begin{equation}
\zeta_A':=\left[s_{\mu_A}^{-1}dx_1\wedge dx_2\wedge dx_3\right].
\end{equation}
\end{defn}
\begin{defn}
For any $k\in \ZZ$, define an $\O_{\overline{M}}[[u]]$-free module $\overline{\H}_{F_A}^{(-k)}$ of rank $\mu_A$ 
as 
\begin{equation}
\overline{\H}_{F_A}^{(-k)}:=
\O_{\overline{M}}[[u]]\cdot u^k\zeta_A'
\bigoplus_{\substack{
1\le i\le 3, \\
1\le j\le a_{i}-1}} 
\O_{\overline{M}}[[u]]\cdot u^kx_i^j\zeta_A'
\bigoplus \O_{\overline{M}}[[u]]\cdot u^ks_{\mu_A}^{-1}x_1x_2x_3\zeta_A',
\end{equation}
where we regard $\zeta_A', x_i^j\zeta_A', s_{\mu_A}^{-1}x_1x_2x_3\zeta_A'$ as 
elements of $\Gamma(\overline{M},\iota_*\H_{F_A}^{(0)})$.
\end{defn}
\begin{prop}\label{extend de Rham}
We have the isomorphism of $\O_M[[u]]$-modules
\begin{equation}
\iota^*\left(\overline{\H}_{F_A}^{(-k)}\right)\simeq\H_{F_A}^{(-k)},\quad i\in\ZZ.
\end{equation}
Namely, $\overline{\H}_{F_A}^{(-k)}$ is an extension of $\H_{F_A}^{(-k)}$ as a free $\O_{\overline{M}}[[u]]$-module.
\end{prop}
\begin{pf}
It is clear by Proposition~\ref{b-lat-k}.
\qed
\end{pf}
\begin{defn}\label{defn-extend-omega F}
Define an $\O_{\overline{M}}$-free module $\overline{\Omega}_{F_A}$ of rank $\mu_A$ 
as 
\begin{equation}
\overline{\Omega}_{F_A}:=
\O_{\overline{M}}\cdot r^{(0)}(\zeta_A')
\bigoplus_{\substack{
1\le i\le 3, \\
1\le j\le a_{i}-1}} 
\O_{\overline{M}}\cdot r^{(0)}(x_i^j\zeta_A')
\bigoplus \O_{\overline{M}}\cdot r^{(0)}(s_{\mu_A}^{-1}x_1x_2x_3\zeta_A'),
\end{equation}
where we regard $r^{(0)}(\zeta_A'), r^{(0)}(x_i^j\zeta_A'), r^{(0)}(s_{\mu_A}^{-1}x_1x_2x_3\zeta_A')$ as elements 
of $\Gamma(\overline{M},\iota_*\Omega_{F_A})$.
\end{defn}
\begin{prop}\label{extend omega F}
We have the isomorphism of $\O_M$-modules
\begin{equation}
\iota^*\left(\overline{\Omega}_{F_A}\right)\simeq\Omega_{F_A}.
\end{equation}
Namely, $\overline{\Omega}_{F_A}$ is an extension of $\Omega_{F_A}$ as a free $\O_{\overline{M}}$-module.
\end{prop}
\begin{pf}
It is almost clear.
\qed
\end{pf}
\begin{prop}\label{rank1-extend-omega F}
$\overline{\Omega}_{F_A}$ is a free ${p}_{*}\O_{\overline{\C}}$-module of rank one.
\end{prop}
\begin{pf}
A ${p}_*\O_{\overline{\C}}$-free base of $\overline{\Omega}_{F_A}$ can be chosen as $r^{(0)}(\zeta_A')$ in 
Definition \ref{defn-extend-omega F}.
\qed
\end{pf}
\begin{prop}\label{extend-b-lat}
We have the following short exact sequence of $\O_{\overline{M}}$-modules
\begin{equation}
0 \rightarrow \overline{\H}^{(-1)}_{F_A}\hookrightarrow \overline{\H}^{(0)}_{F_A} \xrightarrow{\overline{r}^{(0)}} \overline{\Omega}_{F_A}\rightarrow 0.
\end{equation}
\end{prop}
\begin{pf}
It is clear by their definitions.
\qed
\end{pf}
\subsection{Gau\ss--Manin connection}
We define the free $\O_M[[u]]$-module $\T_{\CC_u\times M}$ of rank $\mu_A+1$ as follows:
\begin{equation}
\T_{\CC_u\times M}:=\O_M[[u]]\frac{d}{du}\oplus \O_M[[u]]\otimes_{\O_M}\T_M.
\end{equation}
\begin{defn}
We define a connection, called the {\it Gau\ss--Manin connection},
\begin{equation}
\nabla: \T_{\CC_u\times M} \otimes_{\O_M} \H_{F_A} \rightarrow \H_{F_A}
\end{equation}
by letting; for $\delta \in \T_{M}$ and $\zeta =[\phi dx_1\wedge \dots \wedge dx_n]\in \H_{F_A}$,
\begin{subequations}
\begin{align}
\nabla_{\delta} \zeta &:=[(\frac{1}{u}(\delta F_A)\phi +\delta (\phi ) )dx_1\wedge \dots \wedge dx_n],\\
\nabla_{\frac{d}{du}} \zeta &:=[(-\frac{1}{u^{2}} F_A\phi +\frac{d\phi}{du})dx_1\wedge \dots \wedge dx_n].
\end{align}
\end{subequations}
\end{defn}
\begin{prop}\label{Gauss--Manin connection}
{\it Gau\ss--Manin connection}
$\nabla:\T_{\CC_u\times M}\otimes_{\O_M}\H_{F_A}\longrightarrow \H_{F_A}$ satisfies following;
\begin{enumerate}
\item $\nabla$ is integrable$:$
\begin{equation}
\left[\nabla_{\frac{d}{du}},\nabla_{\frac{d}{du}}\right]=0,\ 
\left[\nabla_{\frac{d}{du}}, \nabla_\delta \right]=0,\ 
\left[\nabla_\delta, \nabla_{\delta'}\right]=
\nabla_{[\delta,\delta']},\ \delta,\delta'\in\T_M.
\end{equation}
\item $\nabla$ satisfies {\it Griffith transversality}$:$ that is, 
\begin{equation}
\nabla:\T_M\otimes_{\O_M}\H_{F_A}^{(-k)}\longrightarrow\H_{F_A}^{(-k+1)},\quad k\in\ZZ.
\end{equation}
\item The covariant differentiation $\nabla_{\frac{d}{du}}$ satisfies 
\begin{equation}\label{eq:GM-2}
\nabla_{u\frac{d}{du}}(\H_{F_A}^{(-k)})\subset\H_{F_A}^{(-k+1)},\quad k\in\ZZ.
\end{equation}
\end{enumerate}
\end{prop}
\begin{pf}
See Proposition 4.5 of \cite{st:1} and reference there in.
\qed
\end{pf}
We shall consider the extension of the Gau\ss--Manin connection $\nabla$ on $\H_{F_A}$
to the one on $\overline{\H}_{F_A}$.
Define the free $\O_{\overline{M}}[[u]]$-module $\T_{\CC_u\times \overline{M}}(-\log M_\infty)$ of rank $\mu_A+1$ as follows:
\begin{equation}
\T_{\CC_u\times \overline{M}}(-\log M_\infty):=\O_{\overline{M}}[[u]]\frac{d}{du}\oplus 
\O_{\overline{M}}[[u]]\otimes_{\O_{\overline{M}}}\T_{\overline{M}}(-\log M_\infty).
\end{equation}
\begin{prop}\label{extension GM connection}
The Gau\ss--Manin connection $\nabla$ on $\H_{F_A}$ extends to $\overline{\H}_{F_A}$ 
with logarithmic poles along $M_\infty$, namely, $\nabla$ induces the connection
\begin{equation}
\overline{\nabla}:\T_{\CC_u\times \overline{M}}(-\log M_\infty)\otimes_{\O_{\overline{M}}}\overline{\H}_{F_A}\longrightarrow
\overline{\H}_{F_A},
\end{equation}
satisfying the following conditions$:$
\begin{enumerate}
\item $\overline{\nabla}$ is integrable$:$
\begin{equation}
\left[\overline{\nabla}_{\frac{d}{du}},\overline{\nabla}_{\frac{d}{du}}\right]=0,\ 
\left[\overline{\nabla}_{\frac{d}{du}}, \overline{\nabla}_\delta \right]=0,\ 
\left[\overline{\nabla}_\delta, \overline{\nabla}_{\delta'}\right]=
\overline{\nabla}_{[\delta,\delta']},\ \delta,\delta'\in\T_{\overline{M}}(-\log M_\infty).
\end{equation}
\item $\overline{\nabla}$ satisfies {\it Griffith transversality}$:$ that is, 
\begin{equation}
\overline{\nabla}:\T_{\overline{M}}(-\log M_\infty)\otimes_{\O_{\overline{M}}}
\overline{\H}_{F_A}^{(-k)}\longrightarrow\overline{\H}_{F_A}^{(-k+1)},\quad k\in\ZZ.
\end{equation}
\item The covariant differentiation $\overline{\nabla}_{\frac{d}{du}}$ satisfies 
\begin{equation}
\overline{\nabla}_{u\frac{d}{du}}(\overline{\H}_{F_A}^{(-k)})\subset\overline{\H}_{F_A}^{(-k+1)},\quad k\in\ZZ.
\end{equation}
\end{enumerate}
\end{prop}
\begin{pf}
We shall check that 
\begin{subequations}
\begin{equation}\label{extend Gauss manin mu}
s_{\mu_A}\nabla_{\frac{\p}{\p s_{\mu_A}}}\left(\overline{\H}_{F_A}\right)\subset \overline{\H}_{F_A}.
\end{equation}
\begin{equation}\label{extend Gauss manin s}
\nabla_{\frac{\p}{\p s_{1}}}\left(\overline{\H}_{F_A}\right)\subset \overline{\H}_{F_A}, \quad 
\nabla_{\frac{\p}{\p s_{i,j}}}\left(\overline{\H}_{F_A}\right)\subset \overline{\H}_{F_A}, \quad
i=1,2,3, \ j=1,\dots, a_{i}-1.
\end{equation}
\begin{equation}\label{extend Gauss manin u}
u\nabla_{\frac{d}{d u}}\left(\overline{\H}_{F_A}\right)\subset \overline{\H}_{F_A}.
\end{equation}
\end{subequations}
First we shall check the condition \eqref{extend Gauss manin mu}. One has
\begin{subequations}
\begin{equation*}\label{eGM mu1}
s_{\mu_A}\nabla_{\frac{\p}{\p s_{\mu_A}}}\zeta_A'=\frac{1}{u}s^{-1}_{\mu_A}x_1x_2x_3\zeta_A'-\zeta_A',
\end{equation*}
\begin{equation*}\label{eGM mu2}
s_{\mu_A}\nabla_{\frac{\p}{\p s_{\mu_A}}}x_{i}^{j}\zeta_A'=\frac{1}{u}s^{-1}_{\mu_A}x_1x_2x_3\cdot x_{i}^{j}\zeta_A'-x_{i}^{j}\zeta_A',
\end{equation*}
\begin{equation*}\label{eGM mu3}
s_{\mu_A}\nabla_{\frac{\p}{\p s_{\mu_A}}}s^{-1}_{\mu_A}x_1x_2x_3\zeta_A'=\frac{1}{u}(s^{-1}_{\mu_A}x_1x_2x_3)^{2}\zeta_A'-2s_{\mu_A}^{-1}x_1x_2x_3\zeta_A'.
\end{equation*}
\end{subequations}
The images of $s^{-1}_{\mu_A}x_1x_2x_3\cdot x_{i}^{j}\zeta_A'$ and $(s^{-1}_{\mu_A}x_1x_2x_3)^{2}\zeta_A'$ by $r^{(0)}$ can be extended to 
$\overline{\Omega}_{F_A}$.
By Proposition~\ref{rank1-extend-omega F} and 
Proposition~\ref{extend-b-lat}, we can show that 
$1/u\cdot s^{-1}_{\mu_A}x_1x_2x_3\zeta_A', \ 1/u\cdot (s^{-1}_{\mu_A}x_1x_2x_3)^{2}\zeta_A' \in \overline{\H}_{F_A}$.    
Next we shall check that the condition \eqref{extend Gauss manin u}. 
The assertion for $s_1$ is obvious. One has
\begin{equation*}
\nabla_{\frac{\p}{\p s_{i',j'}}} \phi \cdot \zeta_A'=\frac{1}{u} \ x_{i'}^{j'}\cdot \phi \cdot \zeta_A',
\end{equation*}
where $\phi$ is an element of the set $\left\{1, \p F_A/\p s_{i,j} \ (i=1,2,3. \ j=1,\dots, a_{i}-1), \ \p F_A/\p s_{\mu_A}\right\}$.
The image of $x_{i'}^{j'}\cdot \phi \cdot \zeta_A$ by $r^{(0)}$ can be extended to 
$\overline{\Omega}_{F_A}$. By Proposition~\ref{rank1-extend-omega F} and 
Proposition~\ref{extend-b-lat}, we can show that 
$1/u\cdot \ x_{i'}^{j'}\cdot \phi \cdot \zeta_A' \in \overline{\H}_{F_A}$.
Finally we shall check the condition \eqref{extend Gauss manin u}. One has
\begin{equation*}
u\nabla_{\frac{d}{d u}} \phi \cdot \zeta_A'=-\frac{1}{u}F_A \cdot \phi \cdot \zeta_A',
\end{equation*}
where $\phi$ is an element of the set $\left\{1, \p F_A/\p s_{i,j} \ (i=1,2,3. \ j=1,\dots, a_{i}-1), \ \p F_A/\p s_{\mu_A}\right\}$.
The image of $F_A \cdot \phi \cdot \zeta_A'$ by $r^{(0)}$ can be extended to 
$\overline{\Omega}_{F_A}$ since we have 
\begin{equation*}
F_A=EF_A+\sum_{i=1}^3\frac{1}{a_i}x_i\frac{\p F_A}{\p x_i}.
\end{equation*}
By Proposition~\ref{rank1-extend-omega F} and 
Proposition~\ref{extend-b-lat}, we can show that $-1/u\cdot F_A \cdot \phi \cdot \zeta_A' \in \overline{\H}_{F_A}$.
The conditions {\rm (i)}, {\rm (ii)} and {\rm (iii)} follow from Proposition~\ref{Gauss--Manin connection}.
\qed
\end{pf}
\subsection{Higher residue pairing}

\begin{defn}
Define an $\O_M$-bilinear form $J_A$ on $\Omega_{F_A}$ by
\begin{equation} 
J_{F_A}(\omega_1, \omega_2):=
- {\rm Res}_{\X/M}\left[
\begin{gathered}
\phi_1\phi_2dx_1\wedge dx_2 \wedge dx_3\\
\frac{\p F_A}{\p x_1}\ \frac{\p F_A}{\p x_2}\ \frac{\p F_A}{\p x_3}
\end{gathered}
\right],
\end{equation}
where $\omega_1 = [\phi_1dx_1\wedge dx_2 \wedge dx_3]$ and $\omega_2 = [\phi_2dx_1\wedge dx_2 \wedge dx_3]$.
\end{defn}
\begin{lem}\label{Gro}
The $\O_M$-bilinear form $J_{F_A}$ on $\Omega_{F_A}$ is non-degenerate.
\end{lem}
\begin{pf}
This is a well-known fact (cf. Section 10.4 of \cite{her:1}). 
\qed
\end{pf}
In order to define the higher residue pairing, we prepare a notation. 
For $P=\sum_{i\in\ZZ} p_iu^i\in \O_M((u))$, 
set $P^*:=\sum_{i\in\ZZ} p_i(-u)^i$ such that $(P^*)^*=P$.
\begin{defn}
An $\O_M$-bilinear form
\begin{equation}
K_{F_A}:\H_{F_A}\otimes_{\O_M}\H_{F_A}\to\O_M((u))
\end{equation}
is called the {\it higher residue pairing} if it satisfies the following properties$:$
\begin{enumerate}
\item For all $\omega_1,\omega_2\in\H_{F_A}$,  
\begin{equation}\tag{K1}\label{K1}
K_{F_A}(\omega_1,\omega_2)= (-1)^3K_{F_A}(\omega_2,\omega_1)^*.
\end{equation}
\item For all $P\in\O_M((u))$ and $\omega_1,\omega_2\in\H_F$, 
\begin{equation}\tag{K2}\label{K2}
PK_{F_A}(\omega_1,\omega_2)=K_{F_A}(P\omega_1,\omega_2)=K_{F_A}(\omega_1,P^*\omega_2).
\end{equation}
\item For all $\omega_1, \omega_2\in\H_{F_A}^{(0)}$,  
\begin{equation}\tag{K3}\label{K3}
K_{F_A}(\omega_1,\omega_2)\in u^{3}\O_M[[u]].
\end{equation}
\item The following diagram is commutative$:$
\begin{equation*}
\begin{CD}
K_{F_A}: & \H_{F_A}^{(0)}\times\H_{F_A}^{(0)} & @>>> & u^3\O_M[[u]]\\
& @VVV & & @VV
{{\rm mod}~u^{4}\O_M[[u]]}V\\
J_{F_A}: & \Omega_{F_A}\times\Omega_{F_A} & @>>> & u^3\O_M.
\end{CD}
\end{equation*}
\item For all $\omega_1, \omega_2\in\H_F$ and $\delta\in\T_M$, 
\begin{equation}\tag{K4}\label{K4}
\delta K_{F_A}(\omega_1,\omega_2)=K_{F_A}(\nabla_\delta\omega_1,\omega_2)+
K_{F_A}(\omega_1,\nabla_\delta\omega_2).
\end{equation}
\item For all $\omega_1, \omega_2\in\H_F$,
\begin{equation}\tag{K5}\label{K5}
u\frac{d}{du}K_{F_A}(\omega_1,\omega_2)=K_{F_A}(u\nabla_\frac{d}{du}\omega_1,\omega_2)
+K_{F_A}(\omega_1,u\nabla_\frac{d}{du}\omega_2).
\end{equation}
\end{enumerate}
\end{defn}
\begin{defn}
Define $K^{(k)}_{F_A}$ for $k\in\ZZ$ by the coefficient of the expansion of $K_{F_A}$ in $u$
\begin{equation}
K_{F_A}(\omega_1,\omega_2):=
\sum_{k\in\ZZ}K^{(k)}_{F_A}(\omega_1,\omega_2)u^{k+3},
\end{equation}
and call it the $k$-th {\it higher residue pairing}.
\end{defn}
\begin{thm}[K.~Saito \cite{Saito}]\label{exsist-higher-residue}
There exists a unique higher residue pairing $K_{F_A}$.
\end{thm}
\begin{pf}
See Theorem of \cite{Saito}.
\qed
\end{pf}
We shall consider the extension $\overline{K}_{F_A}$ of $K_{F_A}$ on $\overline{\H}_{F_A}$.
\begin{prop}\label{extension higher}
The pairing $K_{F_A}$ on $\H_{F_A}^{(0)}$ induces an $\O_{\overline{M}}$-bilinear form 
\begin{equation}
\overline{K}_{F_A}:\overline{\H}_{F_A}^{(0)}\otimes_{\O_{\overline{M}}}\overline{\H}_{F_A}^{(0)}
\longrightarrow u^3\O_{\overline{M}}[[u]]
\end{equation}
whose restriction to $\iota^*\left(\overline{\H}_{F_A}^{(0)}\right)\simeq\H_{F_A}^{(0)}$ 
coincides with $K_{F_A}$.
\end{prop}
\begin{pf}
This follows from Lemma~3.4 of \cite{her:2}, where 
$\overline{M}$, $M$, $M_\infty$, $\H_{F_A}^{(0)}$, $K_{F_A}$, $\overline{\H}_{F_A}^{(0)}$ 
correspond to $X$, $Y$, $D$, $\H$, $P$, ${}_{1}\F$ in \cite{her:2}, respectively.
\qed
\end{pf}
\subsection{Primitive form}
\begin{defn}\label{definition primitive form}
An element $\zeta\in\Gamma(M,\H_{F_A}^{(0)})$ is called a {\it primitive form} for the tuple $(\H_{F_A}^{(0)},\nabla, K_{F_A})$ 
if it satisfies following five conditions;
\begin{enumerate}
\item $u\nabla_{e}\zeta=\zeta$ and $\zeta$ induces $\O_M$-isomorphism:
\begin{equation}\tag{P1}\label{P1}
\T_M[[u]]\simeq \H_{F_A}^{(0)},\quad \sum_{k=0}^\infty\delta_k u^{k}\mapsto 
\sum_{k=0}^\infty u^{k}(u\nabla_{\delta_k}\zeta).
\end{equation}
\item For all $\delta,\delta'\in\T_M$,  
\begin{equation}\tag{P2}\label{P2}
K_{F_A}(u\nabla_\delta\zeta,u\nabla_{\delta'}\zeta)\in \CC\cdot u^{3}.
\end{equation}
\item There exists $r\in\CC$ such that
\begin{equation}\tag{P3}\label{P3}
\nabla_{u\frac{d}{du}+E}\zeta=r\zeta.
\end{equation}
\item There exists a connection $\ns:\T_M\times \T_M\longrightarrow \T_M$ such that 
\begin{equation}\tag{P4}\label{P4}
u\nabla_\delta\nabla_{\delta'}\zeta=\nabla_{\delta\circ\delta'}\zeta
+u\nabla_{\ns_\delta\delta'}\zeta,\quad \delta,\delta'\in\T_M.
\end{equation}
\item There exists an $\O_M$-endomorphism $N:\T_M\longrightarrow \T_M$ such that 
\begin{equation}\tag{P5}\label{P5}
u\nabla_\frac{d}{du}(u\nabla_\delta\zeta)=-\nabla_{E\circ \delta}\zeta+u\nabla_{N\delta}\zeta,\quad 
\delta\in\T_M.
\end{equation}
\end{enumerate}
In particular, the constant $r$ of {\rm (P3)} is called the {\it minimal exponent}.
\end{defn}

It was claimed by K. Saito \cite{S1202-Saito} and is proven by M.~Saito \cite{m.saito} that, for the universal unfolding of 
a cusp singularity, there exists a unique primitive form up to a constant factor with the minimul exponent $r=1$, whose associated exponents are given by the mixed Hodge structure of the cusp singularity. 
In order to prove the mirror isomorphism, we give a its local expression as follows: 
\begin{thm}[cf. \cite{S1202-Saito, m.saito}]\label{primitive form}
Assume that $\chi_A<0$.
There exists a unique primitive form $\zeta_A$ for the tuple $(\H_{F_A}^{(0)},\nabla, K_{F_A})$ 
with the minimal exponent $r=1$ such that 
\begin{equation}
\zeta_A|_{{\bf s}={\bf 0}}=[s_{\mu_A}^{-1}dx_1 \wedge dx_2 \wedge dx_3].
\end{equation}
\end{thm}
\section{Proof of Theorem~\ref{primitive form}}

In this section, we shall prove Theorem \ref{primitive form}. The detail of the proof is not
necessary for following sections. The reader interested in the mirror isomorphism and the period mapping of the primitive form
can skip this section.
 
Set $\overline{M^{s}}:=\{({\bf s},s_{\mu_A})\in \overline{M}~|~{\bf s}={\bf 0}\}$ and
define elements of $\overline{\H}_{F_A^s}^{(0)}:=\left.\overline{\H}_{F_A}^{(0)}\right|_{\overline{M^{s}}}$ by
\begin{subequations}
\begin{align}
\zeta_1:=&[s_{\mu_A}^{-1}dx_1 \wedge dx_2 \wedge dx_3], \\
\zeta_{i,j}:=&[s_{\mu_A}^{-1}x_i^j dx_1 \wedge dx_2 \wedge dx_3],\ i=1,2,3,\ j=1,\dots, a_i-1, \\ 
\zeta_{\mu_A}:=&u\nabla_{s_{\mu_A}\frac{\p}{\p s_{\mu_A}}}\zeta_1,
\end{align}
\end{subequations}
and put 
\begin{equation}
V:=\CC\zeta_1
\bigoplus_{\substack{
1\le i\le 3, \\
1\le j\le a_{i}-1}} 
\CC\zeta_{i,j}
\bigoplus  \CC\zeta_{\mu_A}.
\end{equation}
It is obvious that there is an $\O_{\overline{M^{s}}}$-isomorphism
\begin{equation}
\overline{\H}_{F_A^s}^{(0)}\simeq \O_{\overline{M^s}}[[u]]\otimes_\CC V.
\end{equation}

From now on, we shall check that the elemnts above define
a good section in the sense of Kyoji Saito (see Section 4 of \cite{S1202-Saito}):
\begin{prop}
The elements
$\zeta_1, \ \zeta_{i,j} \ (i=1,2,3, \ j=1,\dots a_{i}-1), \ \zeta_{\mu_A}\in
\overline{\H}_{F_A^s}^{(0)}$
satisfy the following conditions:
\begin{enumerate}
\item $\left.K_{F_A^s}\left(\zeta_i ,\zeta_j\right)\right|_{s_{\mu_A}=0}\in \CC u^3$,
\item $\left.\left(\nabla_{u\frac{d}{du}} \vec{\zeta}\right)\right|_{s_{\mu_A=0}}
=N_{0}\left(\left.u^{-1}\cdot\vec{\zeta}\right|_{s_{\mu_A=0}}\right)+S_{0}\left(\left.\vec{\zeta}\right|_{s_{\mu_A=0}}\right)$,
\end{enumerate}
where $\vec{\zeta}:={}^t (\zeta_1, \zeta_{1,1}, \dots, \zeta_{3,a_{3}-1}, \zeta_{\mu_A})$,
$N_{0}\in M(\mu_A,\QQ)$ and $S_{0}\in M(\mu_A,\CC)$ with $N_{0}$ nilpotent and
$S_{0}$ diagonal. Moreover, $S_{0}$ has eigenvalues 
$\{\alpha_{1},\alpha_{1,1},\dots,\alpha_{3,a_{3}-1},\alpha_{\mu_A}\}$
which coincide with the exponents of the cusp polynomial $f_A$:
\begin{equation}
\alpha_{1}=1, \
\alpha_{i,j}=1+\frac{j}{a_{i}} \ (i=1,2,3, \ j=1,\dots, a_{i}-1), \
\alpha_{\mu_A}=2.
\end{equation}

\end{prop}
\begin{lem}\label{lem:degree of sections}
We have
\begin{subequations}
\begin{align}
&\nabla_{u\frac{d}{du}+\chi_As_{\mu_A}\frac{\p}{\p s_{\mu_A}}}\zeta_1=1\cdot \zeta_1, \\
&\nabla_{u\frac{d}{du}+\chi_As_{\mu_A}\frac{\p}{\p s_{\mu_A}}}\zeta_{i,j}
=\left(1+\frac{j}{a_i}\right)\cdot \zeta_{i,j}, \quad i=1,2,3,\ j=1,\dots, a_i-1,\\
&\nabla_{u\frac{d}{du}+\chi_As_{\mu_A}\frac{\p}{\p s_{\mu_A}}}\zeta_{\mu_A}
=2\cdot \zeta_{\mu_A}.
\end{align}
\end{subequations}
\end{lem}
\begin{pf}
They easily follow from Lemma 3.4 of \cite{ist:2}.
\qed
\end{pf}

\begin{lem}\label{lem:deg of hr}
We have 
\begin{equation}
\deg \left(K^{(k)}_{F_A^s}\left(\zeta_i ,\zeta_j\right)\right)
=\deg(\zeta_i)+\deg(\zeta_j)-3-k,\quad k\in\ZZ.
\end{equation}
\end{lem}
\begin{pf}
The equations \eqref{K4} and \eqref{K5} yields
\begin{align*}
&\ \left(u\frac{d}{du}+\chi_As_{\mu_A}\frac{\p}{\p s_{\mu_A}}\right)
K_{F_A^s}\left(\zeta_i ,\zeta_j\right)\\
=&\ K_{F_A^s}\left(\nabla_{u\frac{d}{du}+\chi_As_{\mu_A}\frac{\p}{\p s_{\mu_A}}}\zeta_i ,\zeta_j\right)
+K_{F_A^s}\left(\zeta_i ,\nabla_{u\frac{d}{du}+\chi_As_{\mu_A}\frac{\p}{\p s_{\mu_A}}}\zeta_j\right).
\end{align*}
The statement now easily follows from the definition of $K^{(k)}_{F_A}$ and Lemma~\ref{lem:degree of sections}. 
\qed
\end{pf}

\begin{lem}\label{lem: higher vanish}
We have 
\begin{equation}
K_{F_A^s}\left(\zeta_i ,\zeta_j\right)\in u^3\O_{\overline{M}^s}[[u]].
\end{equation}
In particular, we have 
\begin{equation}
\left.K_{F_A^s}\left(\zeta_i ,\zeta_j\right)\right|_{s_{\mu_A}=0}\in \CC u^3.
\end{equation}
\end{lem}
\begin{pf}
By Lemma~\ref{lem:deg of hr}, we have 
$\deg \left(K^{(k)}_{F_A^s}\left(\zeta_i ,\zeta_j\right)\right)<0$ if $k\ge 1$ except for the case 
$\zeta_i=\zeta_j=\zeta_{\mu_A}$ and $k=1$.
However, $K^{(1)}_{F_A^s}\left(\zeta_{\mu_A} ,\zeta_{\mu_A}\right)=0$ since $K^{(1)}_{F_A^s}$ 
is skew-symmetric.
Note that $\deg(s_{\mu_A})=\chi_A<0$. 
Therefore, if $k\ge 1$, then $K^{(k)}_{F_A^s}\left(\zeta_i ,\zeta_j\right)\in \O_{\overline{M}^s}$, which 
vanishes at $s_{\mu_A}=0$.
\qed
\end{pf}

\begin{lem}\label{nilpotent}
The residue endomorphism $\left[\nabla_{s_{\mu_A}\frac{\p}{\p s_{\mu_A}}}\right]$ 
on $\left.\overline{\H}_{F_A^{s}}^{(1)}\right|_{s_{\mu_A}=0}$ has zero eigenvalues
and $2\times 2$ nilpotent Jordan block, namely, 
\begin{equation}
\nabla_{s_{\mu_A}\frac{\p}{\p s_{\mu_A}}}\nabla_{s_{\mu_A}\frac{\p}{\p s_{\mu_A}}}\zeta_1=0.
\end{equation}

\end{lem}
\begin{pf}
This lemma follows from Lemma 10.2 in \cite{her:1}. 
\qed
\end{pf}

Hence, combing Lemma~\ref{lem:degree of sections} and Lemma \ref{lem: higher vanish} with Lemma~\ref{nilpotent},
it turns out that the set $\{\zeta_1,\dots, \zeta_{\mu_A}\}$ of 
elements in $\left.\overline{\H}_{F_A^s}^{(0)}\right|_{s_{\mu_A}=0}$
defines a good section, 
which gives a primitive form $\zeta_A$ for the tuple $(\H_{F_A}^{(0)},\nabla, K_{F_A})$ 
with the minimal exponent $r=1$.
Since the exponents defined by this primitive form $\zeta_A$ is given by 
\begin{equation}
\left\{
1,1+\frac{1}{a_1},\dots, 1+\frac{a_1-1}{a_1},
1+\frac{1}{a_2},\dots, 1+\frac{a_2-1}{a_2},
1+\frac{1}{a_3},\dots, 1+\frac{a_3-1}{a_3},2\right\}
\end{equation}
and they coincide with the exponents defined by the mixed Hodge structure for $f_A$.
Up to a constant factor, 
$\zeta_A$ is the primitive form associated to $f_A$ 
obtained by the general theory developed by M.~Saito (see Theorem in Section 3.6 of \cite{m.saito}).
Therefore, we see that
\begin{equation}
\zeta_A|_{{\bf s}={\bf 0}}=\zeta_1=[s_{\mu_A}^{-1}dx_1 \wedge dx_2 \wedge dx_3].
\end{equation}
This finishes the proof of Theorem~\ref{primitive form}.
\section{Mirror Symmetry}

In this section, we shall show the isomorphism of Frobenius manifolds between 
the one constructed from the orbifold Gromov-Witten theory of $\PP^{1}_A$ for the case $\chi_A< 0$
and the one constructed from the pair of the cusp singularity $f_A$ and the primitive form $\zeta_A$
obtained in Theorem \ref{primitive form}.  
\subsection{Mirror isomorphism}

The main theorem of this section is the following:
\begin{thm}\label{cor:classical mirror}
Assume that $\chi_A<0$.
There exists an isomorphism of Frobenius manifolds between the one constructed 
from the Gromov--Witten theory 
for $\PP^{1}_A$ and the one constructed from the pair $(f_A,\zeta_A)$. 
\end{thm}

\begin{pf}
Theorem \ref{cor:classical mirror} immediately follows from Theorem 3.1 in \cite{ist:1}
and the following Theorem \ref{thm:satisfies IST1}.
\qed
\end{pf}

In our previous paper \cite{ist:1}, it is shown that a Frobenius structure with certain conditions can be reconstructed uniquely
and the one constructed from the Gromov--Witten theory of the orbifold projective line $\PP^{1}_A$ with arbitary triplet of positive integers $A$
satisfies the conditions. We shall see the Frobenius structure constructed from the pair $(f_A,\zeta_A)$ also satisfies the conditions, i.e.,
the following Thereom \ref{thm:satisfies IST1} holds: 
\begin{thm}\label{thm:satisfies IST1}
Assume that $\chi_A<0$.
For the Frobenius structure of rank $\mu_A$ and dimension one constructed from the pair $(f_A,\zeta_A)$, 
there exist flat coordinates $t_1$, $t_{1,1},\dots ,t_{3,a_3-1}$, $t_{\mu_A}$ satisfying the following conditions$:$
\begin{enumerate}
\item 
The unit vector field $e$ and the Euler vector field $E$ are given by
\[
e=\frac{\p}{\p t_1},\ E=t_1\frac{\p}{\p t_1}+\sum_{i=1}^3\sum_{j=1}^{a_i-1}\frac{a_i-j}{a_i}t_{i,j}\frac{\p}{\p t_{i,j}}
+\chi_A\frac{\p}{\p t_{\mu_A}}.
\]
\item 
The non-degenerate symmetric bilinear form $\eta$ on $\T_M$ defined by 
\begin{equation*}
\eta(\delta,\delta'):=K_{F_A}^{(0)}(u\nabla_\delta\zeta_A,u\nabla_{\delta'}\zeta),\quad \delta,\delta'\in\T_M
\end{equation*}
satisfies
\begin{align*}
&\ \eta\left(\frac{\p}{\p t_1}, \frac{\p}{\p t_{\mu_A}}\right)=
\eta\left(\frac{\p}{\p t_{\mu_A}}, \frac{\p}{\p t_1}\right)=1,\\ 
&\ \eta\left(\frac{\p}{\p t_{i_1,j_1}}, \frac{\p}{\p t_{i_2,j_2}}\right)=
\begin{cases}
\frac{1}{a_{i_1}}\quad i_1=i_2\text{ and }j_2=a_{i_1}-j_1,\\
0 \quad \text{otherwise}.
\end{cases}
\end{align*}
\item 
The Frobenius potential $\F_{f_A,\zeta_A}$ satisfies $E\F_{f_A,\zeta_A}|_{t_{1}=0}=2\F_{f_A,\zeta_A}|_{t_{1}=0}$,
\[
\left.\F_{f_A,\zeta_A}\right|_{t_1=0}\in\CC[[t_{1,1}, \dots, t_{1,a_1-1},
t_{2,1}, \dots, t_{2,a_2-1},t_{3,1}, \dots, t_{3,a_3-1},e^{t_{\mu_A}}]].
\]
\item 
The restriction of the Frobenius potential $\F_{f_A,\zeta_A}$ to 
the submanifold $\{t_1=e^{t_{\mu_A}}=0\}$ is given as 
\[
\F_{f_A,\zeta_A}|_{t_1=e^{t_{\mu_A}}=0}=\G^{(1)}+\G^{(2)}+\G^{(3)}, 
\]
where $\G^{(i)}\in \CC[t_{i,1},\dots, t_{i,a_i-1}]$, $i=1,2,3$.
\item 
In the frame $\frac{\p}{\p t_1}$, $\frac{\p}{\p t_{1,1}},\dots, 
\frac{\p}{\p t_{3,a_3-1}}$, $\frac{\p}{\p t_{\mu_A}}$ of $\T_M$,
the product $\circ$ can be extended to the limit $t_1=t_{1,1}=\dots=t_{3,a_3-1}=e^{t_{\mu_A}}=0$.
The $\CC$-algebra obtained in this limit is isomorphic to 
\[
\CC[x_1,x_2,x_3]\left/\left(x_1x_2,x_2x_3,x_3x_1,a_1x_1^{a_1}-a_2x_2^{a_2},
a_2x_2^{a_2}-a_3x_3^{a_3}\right)\right.,
\]
where $\p/\p t_{1,1},\p/\p t_{2,1},\p/\p t_{3,1}$ are mapped to
$x_1,x_2,x_3$, respectively.
\item The term $t_{1,1}t_{2,1}t_{3,1}e^{t_{\mu_A}}$ occurs with the coeffcient one in $\F_{f_A,\zeta_A}$. 
\end{enumerate}
\end{thm}

We shall show Theorem \ref{thm:satisfies IST1} by checking one by one that the Frobenius structure constructed from the pair
$(f_A,\zeta_A)$ satisfies the conditions in following subsections. 

\subsection{Condition {\rm (i)}}

\begin{lem}\label{lem:calculation residue}
We have
\begin{equation}
{\rm Res}_{\X^s/M^s}\left[
\begin{gathered}
1\cdot dx_1\wedge dx_2 \wedge dx_3\\
\frac{\p F_A^s}{\p x_1}\ \frac{\p F_A^s}{\p x_2}\ \frac{\p F_A^s}{\p x_3}
\end{gathered}
\right]=0,
\end{equation}
\begin{equation} 
{\rm Res}_{\X^s/M^s}\left[
\begin{gathered}
x_i^j\cdot dx_1\wedge dx_2 \wedge dx_3\\
\frac{\p F_A^s}{\p x_1}\ \frac{\p F_A^s}{\p x_2}\ \frac{\p F_A^s}{\p x_3}
\end{gathered}
\right]=0,\ i=1,2,3,j=1,\dots, a_i-1,
\end{equation}
and 
\begin{equation}
{\rm Res}_{\X^s/M^s}\left[
\begin{gathered}
x_1x_2x_3\cdot dx_1\wedge dx_2 \wedge dx_3\\
\frac{\p F_A^s}{\p x_1}\ \frac{\p F_A^s}{\p x_2}\ \frac{\p F_A^s}{\p x_3}
\end{gathered}
\right]=s_{\mu_A}^{3}.
\end{equation}
\end{lem}
\begin{pf}
Some elementary calculations of residues yield the statement.
\qed
\end{pf}

\begin{lem}\label{lem:0-th residue}
We have
\begin{equation}
K_{F_A^s}^{(0)}(\zeta_A,u\nabla_{\frac{\p}{\p s_1}}\zeta_A)=0,\ 
K_{F_A^s}^{(0)}(\zeta_A,u\nabla_{\frac{\p}{\p s_{i,j}}}\zeta_A)=0,\ i=1,2,3,j=1,\dots, a_i-1,
\end{equation}
and 
\begin{equation}
K_{F_A^s}^{(0)}(\zeta_A,u\nabla_{s_{\mu_A}\frac{\p}{\p s_{\mu_A}}}\zeta_A)=1.
\end{equation}
\end{lem}
\begin{pf}
The statement follows from Lemma~\ref{lem:calculation residue}.
\qed
\end{pf}
\begin{lem}\label{lem:flat fcn}
The one form $\theta\in\Gamma(M,\Omega_M^1)$ defined by 
\begin{multline}
\theta: =K_{F_A}^{(0)}(\zeta_A,u\nabla_{\frac{\p}{\p s_1}}\zeta_A)ds_1\\
+\sum_{i=1}^3\sum_{j=1}^{a_i-1}K_{F_A}^{(0)}(\zeta_A,u\nabla_{\frac{\p}{\p s_{i,j}}}\zeta_A)ds_{i,j}
+K_{F_A}^{(0)}(\zeta_A,u\nabla_{\frac{\p}{\p s_{\mu_A}}}\zeta_A)ds_{\mu_A},
\end{multline}
is a closed form which is independent from the choice of coordinates on $M$.
Moreover, there exists a flat coordinate $t$ such that $\theta=dt$.
\end{lem}
\begin{pf}
See section 3.3 3) of \cite{S1202-Saito}.
\qed
\end{pf}
Combining these two Lemmas, we have $dt|_{{\bf s}=0}=ds_{\mu_A}/s_{\mu_A}$. 
Therefore, we can choose a flat coordinate $t_{\mu_A}$ such that 
\begin{subequations}
\begin{equation}
e^{t_{\mu_A}}:=s_{\mu_A}\cdot t({\bf s},e^{s_{\mu_A}}),\quad \deg(t({\bf s},e^{s_{\mu_A}}))=0, \ 
t({\bf 0},e^{s_{\mu_A}})=1,
\end{equation}
\begin{equation}
\left.\frac{\p t_{\mu_A}}{\p s_1}\right|_{{\bf s}={\bf 0}}=0,\ 
\left.\frac{\p t_{\mu_A}}{\p s_{i,j}}\right|_{{\bf s}={\bf 0}}=0.
\end{equation}
\end{subequations}
Since $\zeta_1, \dots, \zeta_{\mu_A}$ form 
a $\O_{M^{s}}[[u]]$-basis of ${\H}_{F_A}^{(0)}$ at $({\bf 0},q)\in\overline{M}$,
one can choose other flat coordinates $t_1, t_{1,1}, \dots, t_{3,a_3-1}$ such that
\begin{equation}
\left.t_1\right|_{{\bf s}={\bf 0}}=\left.t_{1,1}\right|_{{\bf s}={\bf 0}}=\dots=
\left.t_{3,a_3-1}\right|_{{\bf s}={\bf 0}}=0
\end{equation}
together with the following normalization$;$
\begin{subequations}\label{eq:s-t}
\begin{align}
&\left.\frac{\p t_1}{\p s_1}\right|_{({\bf s},s_{\mu_A})=({\bf 0},0)}=1,\ 
\left.\frac{\p t_1}{\p s_{i,j}}\right|_{({\bf s},s_{\mu_A})=({\bf 0},0)}=0,\ 
\left.\frac{\p t_1}{\p s_{\mu_A}}\right|_{({\bf s},s_{\mu_A})=({\bf 0},0)}=0,\\ 
&\left.\frac{\p t_{i,j}}{\p s_1}\right|_{({\bf s},s_{\mu_A})=({\bf 0},0)}=0,\ 
\left.\frac{\p t_{i,j}}{\p s_{i',j'}}\right|_{({\bf s},s_{\mu_A})=({\bf 0},0)}=\delta_{ii'}\delta_{jj'},\ 
\left.\frac{\p t_{i,j}}{\p s_{\mu_A}}\right|_{({\bf s},s_{\mu_A})=({\bf 0},0)}=0,
\end{align}
\end{subequations}
where $\delta_{ii'}$ and $\delta_{jj'}$ are Kronecker's deltas.
In particular, flat coordinates $t_1, t_{1,1}, \dots, t_{3,a_3-1}, t_{\mu_A}$ satisfy 
\begin{equation}
e=\frac{\p}{\p t_1},\ E=t_1\frac{\p}{\p t_1}+\sum_{i=1}^3\sum_{j=1}^{a_i-1}\frac{a_i-j}{a_i}t_{i,j}\frac{\p}{\p t_{i,j}}
+\chi_A\frac{\p}{\p t_{\mu_A}},
\end{equation}
which is Condition {\rm (i)}.
\subsection{Condition {\rm (ii)}}
\begin{lem}
We have
\begin{equation}
K_{F_A}^{(0)}(\zeta_A,u\nabla_{\frac{\p}{\p t_1}}\zeta_A)=0,\ 
K_{F_A}^{(0)}(\zeta_A,u\nabla_{\frac{\p}{\p t_{i,j}}}\zeta_A)=0,\ i=1,2,3,j=1,\dots, a_i-1,
\end{equation}
and 
\begin{equation}
K_{F_A}^{(0)}(\zeta_A,u\nabla_{\frac{\p}{\p t_{\mu_A}}}\zeta_A)=1.
\end{equation}
\end{lem}
\begin{pf}
By Lemma~\ref{lem:flat fcn}, we have
\[
dt_{\mu_A} =K_{F_A}^{(0)}(\zeta_A,u\nabla_{\frac{\p}{\p t_1}}\zeta_A)dt_1
+\sum_{i=1}^3\sum_{j=1}^{a_i-1}K_{F_A}^{(0)}(\zeta_A,u\nabla_{\frac{\p}{\p t_{i,j}}}\zeta_A)dt_{i,j}
+K_{F_A}^{(0)}(\zeta_A,u\nabla_{\frac{\p}{\p t_{\mu_A}}}\zeta_A)dt_{\mu_A}.
\]
The statement follows.
\qed
\end{pf}
Note that the pairings to consider are constant since we are dealing with flat coordinates.
Therefore, we can evaluate them along $M^s$.
Moreover, by the normalization~\eqref{eq:s-t}, we have
\begin{equation}
\left.K_{F_A}^{(0)}\left(u\nabla_{\frac{\p}{\p t_{i,j}}}\zeta_A,u\nabla_{\frac{\p}{\p t_{i',j'}}}\zeta_A\right)\right|_{{\bf s}=0}
=\ \left.K_{F_A}^{(0)}\left(u\nabla_{\frac{\p}{\p s_{i,j}}}\zeta_A,u\nabla_{\frac{\p}{\p s_{i',j'}}}\zeta_A\right)
\right|_{{\bf s}=0}.
\end{equation}
The statement follows from the same argument in Subsection 4.2 in \cite{ist:2} 
and Lemma~\ref{lem:limit of Jac}.
\subsection{Condition {\rm (iii)}}

By $dt|_{\overline{M^{s}}}=ds_{\mu_A}/s_{\mu_A}$, we have
the $\O_{\overline{M^{s}}}$-isomorphism:
\begin{equation}
\left. \T_{\overline{M}}\left(-\log M_\infty\right)\right|_{\overline{M^{s}}}
\simeq \O_{\overline{M^{s}}}\frac{\p}{\p t_1}
\bigoplus_{\substack{
1\le i\le 3, \\
1\le j\le a_{i}-1}} 
\O_{\overline{M^{s}}}\frac{\p}{\p t_{i,j}}
\bigoplus \O_{\overline{M^{s}}}\frac{\p}{\p t_{\mu_A}},
\end{equation}
where we restrict the flat coordinates in Condition {\rm (i)} on $\overline{M^{s}}$
and denote them by same characters. 
Then, by Proposition~\ref{prop:KS-extension}, we have the $\O_{\overline{M^{s}}}$-isomorphism:
\begin{equation}\label{KS-restrict}
\left.p_{*}\O_{\overline{C}}\right|_{\overline{M^{s}}}\simeq 
\left. \T_{\overline{M}}\left(-\log M_\infty\right)\right|_{\overline{M^{s}}}.
\end{equation}
Condition {\rm (iii)} follows from \eqref{KS-restrict} and the fact that
we can take the the flat coordinate 
$\left.e^{t_{\mu_A}}\right|_{\overline{M^{s}}}=s_{\mu_A}$. 
\subsection{Condition {\rm (iv)}}
Recall that the ideal in the equation~\eqref{eq:Jacobi-ring extension} restricted to $M_\infty$
is given by
\begin{equation}
\left(x_2x_3,x_3x_1,x_1x_2, H_1({\bf x},{\bf s}),H_2({\bf x},{\bf s})\right), 
\end{equation}
where
\begin{equation}
H_{i}({\bf x}, {\bf s}):=a_{i}x_{i}^{a_i}-a_{i+1}x_{i+1}^{a_{i+1}}+
\sum_{j=1}^{a_i-1}j\cdot s_{i,j}\cdot x_i^j-\sum_{j=1}^{a_{i+1}-1}j\cdot s_{i+1,j}\cdot x_{i+1}^j,\ i=1,2.
\end{equation}
In particular, we have 
\begin{equation}\label{eq:cond iv}
\left.\eta\left(\frac{\p}{\p s_{i,j}},\frac{\p}{\p s_{i',j'}}\right)\right|_{s_{\mu_A}=0}=0,\quad \text{if}\ i\ne i'.
\end{equation}
Note that the connection $\ns$ on $\T_M$ in \eqref{P4} is the Levi--Civita connection with respect to $\eta$
(see Proposition 7.9 and Proposition 7.16 of \cite{st:1}) and that $t_{i,j}$ is the solution of the system 
\begin{subequations}
\begin{align}
&\ns^*dt_{i,j}=0,\ \left.t_{i,j}\right|_{{\bf s}=s_{\mu_A}=0}=0,\\ 
&\left.\frac{\p t_{i,j}}{\p s_1}\right|_{{\bf s}=s_{\mu_A}=0}=0,\ 
\left.\frac{\p t_{i,j}}{\p s_{i',j'}}\right|_{{\bf s}=s_{\mu_A}=0}=\delta_{ii'}\delta_{jj'},\ 
\left.\frac{\p t_{i,j}}{\p s_{\mu_A}}\right|_{{\bf s}=s_{\mu_A}=0}=0,
\end{align}
\end{subequations}
where $\ns^*$ is a connection on $\Omega_M^1$ dual to $\ns$.
Therefore, by \eqref{eq:cond iv}, we have
\begin{equation}\label{eq:normalization 2}
\left.\frac{\p t_{i,j}}{\p s_{i',j'}}\right|_{s_{\mu_A}=0}=0,\quad\text{if}\  i\ne i'.
\end{equation}
The third derivatives of the Frobenius potential with respect to flat coordinates 
are given by residues. For example, we have
\begin{equation}\label{3-point}
\frac{\p^3 \F_{f_A,\zeta_A}}{\p t_{i_1,j_1}\p t_{i_2,j_2}\p t_{i_3,j_3}}=-e^{-2 t_{\mu_A}}{\rm Res}_{\CC^3\times M/M}\left[
\begin{gathered}
\frac{\p F_A}{\p t_{i_1,j_1}}\frac{\p F_A}{\p t_{i_2,j_2}}\frac{\p F_A}{\p t_{i_3,j_3}} dx_1\wedge dx_2 \wedge dx_3\\
\frac{\p F_A}{\p x_1}\ \frac{\p F_A}{\p x_2}\ \frac{\p F_A}{\p x_3}
\end{gathered}
\right].
\end{equation}
Therefore, by using the above description of the ideal, the normalization~\eqref{eq:s-t} and \eqref{eq:normalization 2},
we can show that  
\begin{equation}
\left.\lim_{e^{t_{\mu_A}}\to 0}\frac{\p^n \F_{f_A,\zeta_A}}{\p t_{i_1,j_1}\cdots \p t_{i_n,j_n}}
\right|_{t_1=t_{1,1}=\dots =t_{3,a_3-1}=0} \ne 0\ 
\text{only if}\ i_1=\dots=i_n,
\end{equation}
by induction on $n$.
\subsection{Condition {\rm (v)}}
The condition {\rm (v)} easily follows from the equation~\eqref{eq:Jacobi-ring extension} by 
setting $s_1=s_{1,1}=\dots =s_{3,a_3-1}=s_{\mu_A}=0$ together with the normalization~\eqref{eq:s-t}.
\subsection{Condition {\rm (vi)}}
Note that the coefficient of the term $t_{1,1}t_{2,1}t_{3,1}e^{t_{\mu_A}}$ is given by the limit  
\begin{equation}
\lim_{e^{t_{\mu_A}} \to 0}\left(e^{-t_{\mu_A}}\cdot
\left.\frac{\p^3 \F_{f_A,\zeta_A}}{\p t_{1,1}\p t_{2,1}\p t_{3,1}}\right|_{t_1=t_{1,1}=\dots =t_{3,a_3-1}=0} \right).
\end{equation}
Note also that we have the following formula
\begin{equation}
\frac{\p^3 \F_{f_A,\zeta_A}}{\p t_{1,1}\p t_{2,1}\p t_{3,1}}=-{\rm Res}_{\CC^3\times M/M}\left[
\begin{gathered}
\frac{\p F_A}{\p t_{1,1}}\frac{\p F_A}{\p t_{2,1}}\frac{\p F_A}{\p t_{3,1}} \omega_0({\bf x};{\bf s},{s_{\mu_A}})^{-2}dx_1\wedge dx_2 \wedge dx_3\\
\frac{\p F_A}{\p x_1}\ \frac{\p F_A}{\p x_2}\ \frac{\p F_A}{\p x_3}
\end{gathered}
\right].
\end{equation}
Then, by the normalization~\eqref{eq:s-t}, the above limit is reduced to the calculation of 
the following limit 
\begin{equation}
\lim_{e^{t_{\mu_A}} \to 0}\left(\left.-e^{-t_{\mu_A}}\cdot e^{-2 t_{\mu_A}} {\rm Res}_{\CC^3\times M/M}\left[
\begin{gathered}
x_1x_2x_3 dx_1\wedge dx_2 \wedge dx_3\\
\frac{\p F_A}{\p x_1}\ \frac{\p F_A}{\p x_2}\ \frac{\p F_A}{\p x_3}
\end{gathered}
\right]\right|_{{\bf s}=0}\right),
\end{equation}
which is one.
\section{Periods of primitive forms}

In this section, we assume that $\chi_A\ne 0$.
For simplicity, we shall denote by $M$ the Frobenius manifold $M_{(f_A,\zeta_A)}$ constructed from
 the pair of the cusp polynomial $f_A$ and the primitive form $\zeta_A$,
where $\zeta_A$ is the one given in Theorem \ref{primitive form} for $\chi_A<0$ 
and the one in Theorem 3.1 in \cite{ist:2} for $\chi_A>0$. 
We shall systematically calculate the intersection form of $M_{(f_A, \zeta_A)}$ by the period mappings of the primitive form.

\subsection{Period mappings and Intersection forms}\label{preliminary period}

In this subsection, we shall reformulate the period mappings of primitive forms and intersection forms considered by K. Saito \cite{S1202-Saito}
in our situation. 

Put
\begin{equation}
\check{M}:=
\begin{cases}
\CC\times M \ \ \text{if} \ \ \chi_A>0,\\
\CC_{\delta}\times M \ \ \text{if} \ \ \chi_A<0.
\end{cases}
\end{equation}
We also denote the coordinate of $\check{M}$ by $(w, {\bf s}, s_{\mu_A})$ and set $\delta_{w}:=\p/\p w$. 

\begin{defn}
For any $\kappa \in \CC$, define the $\D_{\check{M}}$-module
\begin{equation}
\M^{(\kappa)}:=\D_{\check{M}}/\I^{(\kappa)}_{\check{M}}
\end{equation}
with relations:
\begin{subequations}
\begin{eqnarray}
\lefteqn{\I^{(\kappa)}_{\check{M}}:=
\sum_{\p, \p' \in \T_M} \D_{\check{M}}P(\p, \p')+\sum_{\p \in \T_M} 
\D_{\check{M}}Q_{\kappa}(\p ),}\\
&P(\p, \p'):=\p \p'-(\p \circ \p') \delta_{w}-\ns_{\p}\p', \quad \p, \p' \in \T_M,\\
&Q_{\kappa}(\p):=(E\circ \p) \delta_{w}-(N-\kappa-1)\p, \quad \p \in \T_M.
\end{eqnarray}
\end{subequations}
where $\ns$ and $N$ are the connection and the $\O_M$-endomorphism defined
in Definition \ref{definition primitive form} respectively.
\end{defn}

\begin{defn}
For a $\D_{\check{M}}$-module $\N$,
we set the $\D_{\check{M}}$-module:
\begin{equation}
Sol(\N):={\rm Hom}_{\D_{\check{M}}}(\N, \O_{\check{M}}).
\end{equation}
\end{defn}

\begin{rem}
Since $\M^{(\kappa)}$ has the generator $1$,
we can identify a solution $g\in Sol(\M^{(\kappa)})$ with $g(1)\in \O_{\check{M}}$.
\end{rem}

\begin{defn}
We set a $\D_{\check{M}}$--module:
\begin{equation}
\displaystyle
\widetilde{\M^{(\kappa)}}:=\D_{\check{M}}/
\left(\I^{(\kappa)}_{\check{M}}+
\D_{\check{M}}\left(E-(1-\kappa)\right)\right).
\end{equation}
\end{defn}
\begin{rem}
Recall here that the minimal exponent $r$ of $\zeta_A$ in Theorem \ref{primitive form} and Theorem 3.1 in \cite{ist:2} is one.
We substitue $r=1$ to the original definition in \cite{S1202-Saito}. 
\end{rem}

Set $\check{\D} :=\varphi(\C)$ where $\C$ is the critical set of $F_A$ and
$\varphi$ is the morphism defined in Section \ref{sec:N and T}. 
Since these are also easily defined for $\chi_A>0$ in a natural way, we omit the detail (see \cite{ist:2}). 
However, note that we have to consider everything globally for the case that $\chi_A>0$
(for example, $\C$ is a submanifold in $\CC^3\times M$).  
\begin{lem}\label{lem:sol}
One has
\begin{align}
\lefteqn{Sol(\M^{(\kappa)})\simeq Sol(\widetilde{\M}^{(\kappa)})\oplus \CC \tau_{\mu_A},
\quad \text{if} \ \kappa=1,}\\
&Sol(\M^{(\kappa)})\simeq Sol(\widetilde{\M}^{(\kappa)})\oplus \CC,
\quad \text{if} \ \kappa\ne 1,
\end{align}
where $\tau_{\mu_A}$ is a function on $\check{M}\backslash\check{\D}$
such that $E \tau_{\mu_A}$ is non-zero constant and $\left. \tau_{\mu_A}\right|_{M\backslash \D}=t_{\mu_A}$.
\end{lem}
\begin{pf}
See Section 5 of \cite{S1202-Saito}.
\qed
\end{pf}

Since $Sol(\M^{(\kappa)})$ always contains the constant function, we have
\begin{equation}
dSol(\M^{(\kappa)}) \simeq Sol(\M^{(\kappa)})/\CC_{\check{M}},
\end{equation}  
where $dSol(\M^{(\kappa)})$ is the image of $Sol(\M^{(\kappa)})$ in 
$\Omega^{1}_{\check{M}}$ under the differential $d$.

\begin{lem}\label{lem:dsol}
There exist an isomorphism:
\begin{equation}
\Omega^{1}_{M, t} \simeq \O_{M, t}\otimes\left(\left. dSol(\M^{(\kappa)})\right|_{w=0}\right)_{t}, \quad {t\in M\backslash \D},
\end{equation}
where $\D$ is the image of $\check{\D}$ by the natural projection to $M$
and coincides with the support set for the kernel of the multiplication by the Euler vector field $E$.
\end{lem}
\begin{pf}
See Section 5 of \cite{S1202-Saito}.
\qed
\end{pf}

Set
\begin{equation}
\X_{(w, {\bf s}, s_{\mu_A})}:= \{(x_1, x_2, x_3)\in \CC^{3}|w-F_A({\bf x}, {\bf s}; s_{\mu_A})=0\}, \ (w, {\bf s}, s_{\mu_A})\in \check{M}.
\end{equation}
We can relate an element $\omega\in \Gamma(M, \H^{(0)}_{F_A})$ with the Gelfand--Leray form:
\begin{equation}
\displaystyle
\check{\omega}:=
{\rm Res}_{\X_{(w, {\bf s}, s_{\mu_A})}}\frac{\left[\omega\right]}{w-F_A({\bf x}, {\bf s}; s_{\mu_A})}\in
\bigcup_{\substack{(w, {\bf s}, s_{\mu_A}) \in \check{M} \backslash \check{\D}}} H^{2}(\X_{(w, {\bf s}, s_{\mu_A})}, \CC),
\end{equation}
where $[\omega]$ is the image of $\omega$ for the
Fourier--Laplace transformation: $u^{-1}\mapsto \delta_{w}$.

Let $\beta(w,{\bf s},s_{\mu_A}) \in \displaystyle\bigcup_{\substack{(w, {\bf s}, s_{\mu_A}) \in \check{M} \backslash \check{\D}}} H_2 (\X_{(w, {\bf s}, s_{\mu_A})}, \ZZ)$ be a horizontal family of homology classes
defined on a simply connected domain of a covering space of $\check{M}\backslash \check{\D}$.
Then, by considering the Gelfand--Leray form of the primitive form $\zeta_A$, 
one can consider the period $\displaystyle\int_{\beta(w,{\bf s},s_{\mu_A})}\check{\zeta_A}$. 

\begin{lem}\label{period mapping}
One has 
\begin{equation}
\displaystyle\int_{\beta(w,{\bf s},s_{\mu_A}) }\check{\zeta_A}\in Sol(\widetilde{\M}^{(1)}).
\end{equation}
\end{lem}

\begin{pf}
Lemma \ref{period mapping} immediately follows from Note 2 of Section 5 in \cite{S1202-Saito}.
\qed
\end{pf}

Here we shall recall the intersection form of a Frobenius manifold defined by Dubrovin \cite{d:1} and
the intersection form considered by K. Saito \cite{S1202-Saito}, and compare them. 
\begin{defn}[cf. \cite{d:1}]\label{intersection of F-manifold}
Let $M$ be a Frobenius manifold and $E$ the Euler vector field of $M$.
For $1$-forms $\omega_1, \omega_2\in \Omega^{1}_M$, we put
\begin{equation}
I_{M}(\omega_1, \omega_2):=i_E(\omega_1\bullet \omega_2),
\end{equation}
where $\bullet$ is the induced the operation of multiplication of tangent vectors on the Frobenius manifold $M$ and the
duality between tangent and cotangent spaces established by the non-degenerate bilinear form $\eta$ of $M$, and
$i_{E}$ is the operator of contraction of a $1$-form with the Euler vector field $E$.
Moreover, by using flat coordinates of $M$, one can reformulate the intersection form $I_{M}$ as follows:
\begin{equation}
I_{M}(dt^{i}, dt^{j})=\sum_{k,l=1}^{\mu}\eta^{ik}\eta^{jl}E(\p_{k}\p_{l}\F_{M}),
\end{equation}
where we denote by $t_{i}$ the flat coordinate of $M$, by 
$\F_{M}$ the Frobenius potential of $M$ and set $\p_{i}:=\p/\p t_{i}$ and $(\eta^{ab}):=(\eta(\p_{a},\p_{b}))^{-1}$.
\end{defn}

\begin{defn}[cf. \cite{S1202-Saito}]\label{intersection form K.S}
We put the $\O_{\check{M}}$-bilinear form
\begin{subequations}
\begin{equation}
I: \Omega^{1}_{\check{M}}\times \Omega^{1}_{\check{M}}\longrightarrow \O_{\check{M}} 
\end{equation}
\begin{equation}
I(\omega_1,\omega_2):=\sum_{a,b=1}^{\mu_A}
i_{\p_{a}}(\omega_1)\cdot
\eta^{ab}\cdot i_{w\p_{b}+E\circ\p_{b}}(\omega_2),
\end{equation}
\end{subequations}
where $i_{\p_a}$ is the operator of contraction of a $1$-form with the vector field $\p_a\in \T_{M}$, and call it the intesection form of $\check{M}$.
\end{defn}
Combining Lemma \ref{lem:dsol} with Definition \ref{intersection of F-manifold} and Definition \ref{intersection form K.S}, 
one can see that the $\O_{\check{M}}$-bilinear form $I$ restricted to $w=0$
coincides the intersection form of Frobenius manifold $I_{M_{(f_A,\zeta_A)}}$.
For the simplicity, we shall denote $I_{M_{(f_A,\zeta_A)}}$ by $I_{(f_A,\zeta_A)}$.

\subsection{Calculations for the intersection form of $M_{(f_A, \zeta_A)}$}

From now on, we denote by $\tilde{\zeta_A}$ the form $\check{\zeta_A}$ restricted to $w=0$.
Namely, we set
\begin{equation}\label{zetatilde}
\displaystyle
\tilde{\zeta_A}:=
{\rm Res}_{\X_{(0, {\bf s}, s_{\mu_A})}}\frac{\left[\zeta_A\right]}{-F_A({\bf x}, {\bf s}; s_{\mu_A})}.
\end{equation}
We shall calculate the intersection form of the Frobenius manifold $I_{(f_A,\zeta_A)}$
by considering the intersections of vanishing cycles in a fiber.
\begin{lem}\label{Fourier-Laplace of primitve form}
The primitive form $\zeta_A$ satisfies the following equation:
\begin{equation}\label{eq:Fourier-Laplace of primitve form}
\frac{1}{(2\pi \sqrt{-1})^{2}}\int_{\gamma_{0}({\bf s}, s_{\mu_A})}\tilde{\zeta_A} =1, 
\end{equation}
where $\gamma_{0}({\bf s}, s_{\mu_A})$ is the horizontal family of the homology in 
$\displaystyle\bigcup_{(0, {\bf s}, s_{\mu_A}) \in \check{M} \backslash \check{\D}} H_{2}(\X_{(0, {\bf s}, s_{\mu_A})},\CC)$
corresponding to the relative $3$-cycle 
$\Gamma_{0}=\{|x_{1}|=|x_{2}|=|x_{3}|=\varepsilon \}\subset \CC^{3}\times M$. 
In particular, if $\chi_A>0$, one has 
\begin{subequations}
\begin{align}
\zeta_{A}=[e^{-t_{\mu_{A}}}dx_{1}\wedge dx_{2}\wedge dx_{3}],\label{acp}\\
\tilde{\zeta_A}={\rm Res}_{F_A=0}\left[
\frac{e^{-t_{\mu_A}} \ dx_{1}\wedge dx_{2}\wedge dx_{3}}{-F_A} \label{gl of acp}
\right].
\end{align}
\end{subequations}
\end{lem}
\begin{pf}
\noindent
\underline{Case {\rm (i)} $\chi_A<0$}

We can evaluate the left hand of the equation \eqref{eq:Fourier-Laplace of primitve form} at ${\bf s}={\bf 0}$.
Then one has
\begin{eqnarray*}
\lefteqn{
\frac{1}{(2\pi \sqrt{-1})^{2}}\int_{\gamma_{0}({\bf s}, s_{\mu_A})}\tilde{\zeta_A} =
\frac{1}{(2\pi \sqrt{-1})^{3}}\int_{\Gamma_{0}}
\frac{s_{\mu_A}^{-1} \ dx_{1}\wedge dx_{2}\wedge dx_{3}}{-f_A}}\\
&&=\frac{1}{(2\pi \sqrt{-1})^{3}}\int_{\Gamma_{0}}
\frac{dx_{1}\wedge dx_{2}\wedge dx_{3}}{x_{1}x_{2}x_{3}}+
\left\{\sum_{n=1}^{\infty}\left(s_{\mu_A}\cdot\frac{\sum^{3}_{i=1}x^{a_{i}}_{i}
}{x_{1}x_{2}x_{3}}\right)^{n}
\right\}\cdot\frac{dx_{1}\wedge dx_{2}\wedge dx_{3}}{x_{1}x_{2}x_{3}}.
\end{eqnarray*}
Here one has
\[
\displaystyle x_1^{a_1\cdot e_{1}}x_2^{a_{2}\cdot e_{2}}x_3^{a_3\cdot e_{3}}\ne (x_1x_2x_3)^{n} \ \ 
\text{if} \ \  e_{i}\ge 0, \ \sum_{i=1}^{3}e_{i}=n
\]
because, if the equality is attained, one has the contradictory inequation:
\[
\displaystyle n>(\sum_{i=1}^{3}\frac{n}{a_i})=\sum_{i=1}^{3} e_{i}= n.
\]
Therefore one has
\begin{equation*}
\frac{1}{(2\pi \sqrt{-1})^{2}}\int_{\gamma_{0}({\bf s}, s_{\mu_A})}\tilde{\zeta_A}=\frac{1}{(2\pi \sqrt{-1})^{3}}\int_{\Gamma_{0}}
\frac{dx_{1}\wedge dx_{2}\wedge dx_{3}}{x_{1}x_{2}x_{3}}=1.
\end{equation*}

\noindent
\underline{Case {\rm (ii)} $\chi_A>0$}

First, we shall show the equations \eqref{acp} and \eqref{gl of acp}.
By Theorem 3.1 in \cite{ist:2}, 
the element $[s_{\mu_A}^{-1}dx_1 \wedge dx_2 \wedge dx_3]\in \H_{F_A}^{(0)}$ is a primitive form. 
By Lemma 4.2 in \cite{ist:2} and Lemma \ref{lem:flat fcn}, we can choose $t_{\mu_A}:=\log s_{\mu_A}$ as a flat coordinate. 
Therefore we have \eqref{acp} and \eqref{gl of acp}.

Next, we shall show the equation \eqref{eq:Fourier-Laplace of primitve form}. By equation \eqref{gl of acp},
one has
\begin{eqnarray*}
\lefteqn{
\frac{1}{(2\pi \sqrt{-1})^{2}}\int_{\gamma_{0}({\bf s}, s_{\mu_A})}\tilde{\zeta_A}=
\frac{1}{(2\pi \sqrt{-1})^{3}}\int_{\Gamma_{0}}
\frac{e^{-t_{\mu_A}} \ dx_{1}\wedge dx_{2}\wedge dx_{3}}{-F_A}}\\
&&=\frac{1}{(2\pi \sqrt{-1})^{3}}\int_{\Gamma_{0}}
\left\{\sum_{n=1}^{\infty}\left(s_{\mu_A}\cdot\frac{\sum^{3}_{i=1}x^{a_{i}}_{i}
+s_1+\sum^{3}_{i=1}\sum^{a_{i}-1}_{j=1}s_{i,j}x_{i}^{j}
}{x_{1}x_{2}x_{3}}\right)^{n}
\right\}\cdot\frac{dx_{1}\wedge dx_{2}\wedge dx_{3}}{x_{1}x_{2}x_{3}}
\\
&&+\frac{1}{(2\pi \sqrt{-1})^{3}}\int_{\Gamma_{0}}
\frac{dx_{1}\wedge dx_{2}\wedge dx_{3}}{x_{1}x_{2}x_{3}}.
\end{eqnarray*}
Here one has
\begin{equation*}
\displaystyle\prod^{a_1}_{j=1}x_1^{j\cdot e_{1,j}}\prod^{a_2}_{j=1}x_2^{j\cdot e_{2,j}}\prod^{a_3}_{j=1}x_3^{j\cdot e_{3,j}}\ne (x_1x_2x_3)^{n}
\quad \text{if} \ \ e_{i,j}\ge 0, \ \ \sum_{i=1}^{3}\sum_{j=1}^{a_{i}}e_{i,j}\le n
\end{equation*}
because, if the equality is attained, one has the following contradictory inequation:
\begin{equation*}
n<(\sum_{i=1}^{3}\frac{n}{a_i})=\sum_{i=1}^{3}\sum_{j=1}^{a_{i}}\frac{j}{a_i}\cdot e_{i,j}\le \sum_{i=1}^{3}\sum_{j=1}^{a_{i}}e_{i,j}\le n.
\end{equation*}
Then one has
\begin{equation*}
\frac{1}{(2\pi \sqrt{-1})^{2}}\int_{\gamma_{0}({\bf s}, s_{\mu_A})}\tilde{\zeta_A}=\frac{1}{(2\pi \sqrt{-1})^{3}}\int_{\Gamma_{0}}
\frac{dx_{1}\wedge dx_{2}\wedge dx_{3}}{x_{1}x_{2}x_{3}}=1.
\end{equation*}
Therefore we have Lemma~\ref{Fourier-Laplace of primitve form}.
\qed
\end{pf}

\begin{lem}\label{isomorphism of lattice}
There exists the isomorphism between lattices$:$
\begin{equation}
(H_{2}(\X_{(w, {\bf s}, s_{\mu_A})},\ZZ), -I_{H_2}) \simeq (\widetilde{\h}^{*}_A, \left<,\right>) \ \text{if} \ (w, {\bf s}, s_{\mu_A})\in \check{M}\backslash \check{\D},
\end{equation}
where $I_{H_2}$ is the intersection form for cycles in the fiber $\X_{(w, {\bf s}, s_{\mu_A})}$.
\end{lem}

\begin{pf}
See \cite{gl} and
Section 3 of \cite{t:2}.
\qed
\end{pf}

Here $T_A$ is the following
Coxeter--Dynkin diagram:
\begin{equation}
\begin{xy}
\xygraph{
    \bullet ([]!{+(0,-.3)} {(1,a_1-1)}) - [r] \cdots - [r]
    \bullet ([]!{+(0,-.3)} {(1,1)}) - [r] 
    \bullet ([]!{+(.3,-.3)} {1}) (
        - [d] \bullet ([]!{+(.5,0)} {(2,1)})( 
        - [d] \vdots - [d] \bullet ([]!{+(.8,0)} {(2,a_2-1)})),
        - [r] \bullet ([]!{+(0,-.3)} {(3,1)})
        - [r] \cdots - [r] \bullet ([]!{+(0,-.3)} {(3,a_3-1)})
)}
\end{xy}
\end{equation}
and $\h_A$ the complexified Cartan subalgebra of the Kac--Moody Lie algebra associated to $T_A$ (in particular, which is the simple Lie algebra for the case that $\chi_A>0$).
Denote by $\alpha_1,\dots, \alpha_{(3, a_3-1)}\in \h_A^{*}:={\rm Hom}_{\CC}(\h_A,\CC)$ simple roots corresponding the vertices in $T_{A}$, 
by $\alpha^{\vee}_1,\dots, \alpha^{\vee}_{(3, a_3-1)}\in \h_A$ simple coroots
and by $\left<, \right>:\h^{*}_A\otimes_{\CC} \h_A\longrightarrow \CC$ the natural pairing.
The Weyl group $W_{A}$ is a group generated by reflections
\begin{equation}
r_{\bf i}(h):=h-\left<\alpha_{\bf i}, h\right>\alpha^{\vee}_{\bf i}, \quad h\in \h_A, \ {\bf i}=1, (1,1),\dots, (i, j),\dots, (3,a_3-1),
\end{equation}
where $\left<, \right>$ denotes the natural pairing
$\left<, \right>:\h^{*}_A\otimes_{\CC} \h_A\longrightarrow \CC$.
Moreover, we set $\widetilde{\h}_A$ is the complexified 
Cartan subalgebra of the affine Lie algebra
associated to $T_A$.

\begin{lem}\label{subspace of ah}
Denote by $\delta\in \widetilde{\h}^{*}_A:={\rm Hom}_{\CC}(\widetilde{\h}_A, \CC)$ the generator of the imaginary root.
Then one has the isomorphism of affine spaces
\begin{equation}
\h_A\simeq \left\{\widetilde{h}\in \widetilde{\h_A} | \left<\widetilde{h}, \delta\right>
=1\right\}
\end{equation}
which is compatible with the action of the affine Weyl group $\widetilde{W}_A$ on both
sides, where the action on the left hand side is defined by
\begin{equation}
h\mapsto w(h)+\sum_{i=1}^{3}\sum^{a_{i}-1}_{j=1}m_{(i,j)}\alpha^{\vee}_{(i,j)},\quad m_{(i,j)} \in \ZZ
\end{equation}
and the one on the right hand side is the natural one.
\end{lem}
\begin{pf}
Some elementary calculations yield the statement.
\qed
\end{pf}

By Lemma \ref{Fourier-Laplace of primitve form}, Lemma \ref{isomorphism of lattice} and Lemma \ref{subspace of ah},
we can identify $\delta$ with the $2$-cycle $\gamma_{0}({\bf s}, s_{\mu_A})$ in the fiber.
We denote by $\gamma_{\bf i}({\bf s},s_{\mu_A})\in \displaystyle\bigcup_{\substack{(0, {\bf s}, s_{\mu_A}) \in \check{M} \backslash \check{\D}}} 
H_2 (\X_{(0, {\bf s}, s_{\mu_A})}, \ZZ)$ 
the image of $\alpha_{\bf i}$ 
by the composition of the isomorphisms in Lemma \ref{isomorphism of lattice} and Lemma \ref{subspace of ah}. 

Under the above notations and the natural identification between $M$ and $\left. \check{M}\right|_{w=0}$, the following is the main theorem in this section:
\begin{thm}\label{calc: intersection form 2}
Consider the periods
\begin{equation}
x_{\bf i}:=\frac{1}{(2\pi\sqrt{-1})^{2}}\int_{\gamma_{\bf i}({\bf s},s_{\mu_A})} \tilde{\zeta_A}, \quad {\bf i}=1, (1,1), \dots, (i,j), \dots, (3, a_{3}-1),
\end{equation}
where $\gamma_{\bf i}({\bf s},s_{\mu_A})$ is the horizontal family of homology classes 
in $\displaystyle\bigcup_{\substack{(0, {\bf s}, s_{\mu_A}) \in \check{M} \backslash \check{\D}}} H_{2}(\X_{(0, {\bf s}, s_{\mu_A})},\CC)$ identified with $\alpha_{\bf i}$, and the function
\begin{equation}
x_{\mu_A}:=\frac{1}{2\pi \sqrt{-1}} t_{\mu_A}=\frac{1}{2\pi \sqrt{-1}} \log s_{\mu_A}.
\end{equation}
They define the flat coordinates with respect to $I_{(f_A,\zeta_A)}$ on the monodromy covering space of $M\backslash \D$.
Moreover, one has
\begin{subequations}
\begin{align}
\lefteqn{I_{{(f_A,\zeta_A)}} (dx_{\bf i}, dx_{\bf j})
=\frac{-1}{(2\pi \sqrt{-1})^2}\left<\alpha_{\bf i},\alpha_{\bf j}^{\vee}\right>,\label{eq:1}}\\
&I_{{(f_A,\zeta_A)}} (dx_{\mu_A}, dx_{\bf i})
=I_{{(f_A,\zeta_A)}} (dx_{\bf i}, dx_{\mu_A})=0,\label{eq:3}\\
&I_{{(f_A,\zeta_A)}}(dx_{\mu_A},dx_{\mu_A})
=\displaystyle \frac{1}{(2\pi \sqrt{-1})^{2}}\chi_A.\label{eq:4}
\end{align}
\end{subequations}
\end{thm}
\begin{pf}
The first assertion immediately follows from Note 2 of Section 5 in \cite{S1202-Saito} and Lemma \ref{lem:sol}.
We shall show the second assertion. 
First, we shall show the equations \eqref{eq:3} and \eqref{eq:4}. 
By Lemma \ref{lem:dsol} and Definition \ref{intersection of F-manifold},  
one has
\begin{subequations}
\begin{align*}
I_{{(f_A,\zeta_A)}}(dt_{\mu_A}, dt_{1})&=
\displaystyle \sum_{\alpha, \beta} \eta^{\mu_A \alpha}
\eta^{1 \beta}E(\p_{\alpha}\p_{\beta}\F_{{(f_A,\zeta_A)}})\\
&=1\cdot 1\cdot E(\p_{1}\p_{\mu_A}\F_{{(f_A,\zeta_A)}})=t_{1},\\
I_{{(f_A,\zeta_A)}}(dt_{\mu_A}, dt_{i,j})&=
\displaystyle \sum_{\alpha, \beta} \eta^{\mu_A \alpha}
\eta^{(i,j) \beta}E(\p_{\alpha}\p_{\beta}\F_{{(f_A,\zeta_A)}})\\
&=1\cdot a_{i}\cdot E(\p_1\p_{(i,a_{i}-j)}\F_{{(f_A,\zeta_A)}})= \frac{a_{i}-j}{a_{i}}t_{i,j},\\
I_{(f_A,\zeta_A)}(dt_{\mu_A},dt_{\mu_A})
&=\displaystyle \sum_{\alpha,\beta} \eta^{\mu_A \alpha}
\eta^{\mu_A \beta}E(\p_{\alpha}\p_{\beta}\F_{(f_A,\zeta_A)})\\
&=1\cdot 1\cdot E(\p_{1}\p_{1}\F_{{(f_A,\zeta_A)}})=\chi_A.
\end{align*}
\end{subequations}
We shall substitute above calculations for following equations. 
By Lemma \ref{period mapping}, one has 
\begin{subequations}
\begin{align*}
\lefteqn{
I_{(f_A,\zeta_A)} (dx_{\mu_A}, dx_{\bf i})
=\displaystyle \sum_{\alpha,
\beta}I_{(f_A,\zeta_A)}(dt_{\alpha},dt_{\beta})\p_{\alpha}x_{\mu_A}\p_{\beta}x_{\bf i}}\\
&&=\sum_{\beta}\frac{1}{2\pi\sqrt{-1}}I_{(f_A,\zeta_A)}
(dt_{\mu_A},dt_{\beta})\p_{\beta}x_{\bf i}
=\frac{1}{2\pi \sqrt{-1}}Ex_{\bf i} 
=0,\\
\lefteqn{
I_{(f_A,\zeta_A)}(dx_{\mu_A},dx_{\mu_A})
=\displaystyle
I_{(f_A,\zeta_A)}\left(\frac{1}{2\pi\sqrt{-1}}dt_{\mu_A},
\frac{1}{2\pi\sqrt{-1}}dt_{\mu_A}\right)
=\displaystyle \frac{1}{(2\pi \sqrt{-1})^{2}}\chi_A.}
\end{align*}
\end{subequations}

Finally we shall show the equation \eqref{eq:1}: 
\[
I_{{(f_A,\zeta_A)}} (dx_{\bf i}, dx_{\bf j})
=\frac{-1}{(2\pi \sqrt{-1})^2}\left<\alpha_{\bf i}, \alpha_{\bf j}^{\vee}\right>.
\]
The equation \eqref{eq:1} immediately follows from Lemma \ref{isomorphism of lattice} and
the following Lemma \ref{calc: intersection form 1}: 
\begin{lem}\label{calc: intersection form 1}
Let $\displaystyle\beta_{\bf i}(w, {\bf s}, s_{\mu_A}) 
\in \bigcup_{(w, {\bf s}, s_{\mu_A}) \in \check{M} \backslash \check{\D}}H_2 (\X_{(w, {\bf s}, s_{\mu_A})}, \ZZ), \ {\bf i}=1,\dots, (3, a_{3}-1)$
be a horizontal family of homology defined on a simply connected domain of a covering space of $\check{M}\backslash \check{\D}$.
Then, one has 
\begin{multline}
\frac{-1}{(2\pi \sqrt{-1})^{2}}\sum_{a,b=1}^{\mu_{A}}
\p_{a}\left(\int_{\beta_{\bf i}(w, {\bf s}, s_{\mu_A})} \check{\zeta_{A}}\right)\cdot
\eta^{ab}\cdot(E\circ\p_{b})\left(\int_{\beta_{\bf j}(w, {\bf s}, s_{\mu_A})} \check{\zeta_A}\right)\\
\displaystyle=I_{H_2}(\beta_{\bf i}(w, {\bf s}, s_{\mu_A}), \beta_{\bf j}(w, {\bf s}, s_{\mu_A})).
\end{multline}
\end{lem}
\begin{pf}
See Theorem 3.4 in \cite{S1202-Saito} (a factor $(-1)$ is missing in the refernce). 
\qed
\end{pf}
Therefore we have Theorem \ref{calc: intersection form 2}.
\qed
\end{pf}
\section{Frobenius manifold $M_{\widehat{W}_A}$}

We shall recall the Frobenius manifold constructed from the invariant theory of 
an extended affine Weyl group by Dubrovin-Zhang in \cite{dz:1}.

Under the assumption $\chi_A\ne 0$, the Cartan matrix for $T_A$ is nondegenarate.
Set $\widehat{\h}_A:=\h_A\times \CC$.
The affinization $\widetilde{W}_A$ of  $W_A$ acts on $\widehat{\h}_A$ by
\begin{equation}
(h, x_{\mu_A})\mapsto (w(h)+\sum_{i=1}^{3}\sum^{a_{i}-1}_{j=1}m_{(i,j)}\alpha^{\vee}_{(i,j)}, x_{\mu_A}),\quad m_{(i,j)} \in \ZZ
\end{equation}
and $\ZZ$ acts on $\widehat{\h}_A$ by 
\begin{equation}
(h, x_{\mu_A})\mapsto (h+m\omega^{\vee}_{1}, x_{\mu_A}+m), \quad m\in \ZZ,
\end{equation}
where 
$\omega^{\vee}_{1}, \omega^{\vee}_{(1,1)}, \dots, \omega^{\vee}_{(3,a_{3}-1)}$
denotes the fundamental coweights, the elements of $\h_A$ sarisfying $\left<\alpha_{\bf i}, \omega^{\vee}_{\bf j}\right>=\delta_{\bf ij}$
(where $\delta_{\bf ij}$ is the Kronecker's delta).
Then, $\widehat{W}_{A}$ is defined as a group acting on $\widehat{\h}_A$
generated by $\widetilde{W}_A$ and $\ZZ$ with the above actions on $\widehat{\h}_A$. In particular,
one has the following exact sequence
\[
1\rightarrow \widetilde{W}\rightarrow \widehat{W}\rightarrow \ZZ \rightarrow 1.
\]

By the invariant theory of $\widehat{W}_A$, Dubrovin–-Zhang \cite{dz:1} give the following:
\begin{thm}[\cite{dz:1}]\label{dz:thm}
Assume that $\chi_A>0$.
There exists a unique {\rm Frobenius} structure of rank $\mu_A$ and dimension one on 
$M_{\widehat{W}_A}:=\widehat{\h}_A/\widehat{W}_A$ with flat coordinates 
$t_{1},t_{1,1},\dots, t_{i,j}, \dots, t_{3,a_{3}-1}, t_{\mu_A}:=(2\pi\sqrt{-1})x_{\mu_A}$ such that
\begin{equation}
e=\frac{\p}{\p t_1}, \ \ E=t_1\frac{\p}{\p t_1}+\sum_{i=1}^3\sum_{j=1}^{a_i-1}\frac{a_i-j}
{a_i}t_{i,j}\frac{\p}{\p t_{i,j}}
+\chi_A\frac{\p}{\p t_{\mu_A}},
\end{equation}
and the intersection form $I_{\widehat{W}_A}$ is given by
\begin{subequations}
\begin{align}
\lefteqn{
I_{\widehat{W}_A}(\alpha_{\bf i}, \alpha_{\bf j})=\frac{-1}{(2\pi\sqrt{-1})^{2}}
\left<\alpha_{\bf i}, \alpha^{\vee}_{\bf j}\right>, \quad {\bf i,j}=1, (1,1),\dots,(3,a_{3}-1),}\\
&I_{\widehat{W}_A}(\alpha_{\bf i},dx_{\mu_A})=I_{\widehat{W}_A}(dx_{\mu_A},\alpha_{\bf i})
=0,\quad {\bf i}=1, (1,1),\dots,(3,a_{3}-1),\\
&I_{\widehat{W}_A}(dx_{\mu_A},dx_{\mu_A})
=\frac{1}{(2\pi\sqrt{-1})^{2}}\chi_A,
\end{align}
\end{subequations}
where we identify the cotangent space of $M_{\widehat{W}_A}$
with $\h^{*}_A\oplus\CC dx_{\mu_A}$.
\end{thm}
\begin{pf}
See Theorem 2.1 in \cite{dz:1}.
\qed
\end{pf}

\section{Isomorphism between $M_{(f_A,\zeta_A)}$ and $M_{\widehat{W}_A}$}

In this section, we shall show the
isomorphism of Frobenius manifolds between
the one constructed from the pair $(f_A,\zeta_A)$ in \cite{ist:2} and
the one constructed from the invariant theory of an extended 
affine Weyl group $\widehat{W}_A$ in \cite{dz:1}.

\subsection{Reconstruction Theorem via Intersection forms} 

The following Theorem \ref{second} might be known to experts. However the complete proof is not found in any literature.
For this reason, we shall give a proof suitable for our situation here.
\begin{thm}\label{second}
A Frobenius manifold $M$ of rank $\mu_A$ and dimension one with the 
following $e$ and $E$
is uniquely determined by the intersection form $I_M$ $:$
\begin{equation}
e=\frac{\p}{\p t_1},\ E=t_1\frac{\p}{\p t_1}+\sum_{i=1}^3\sum_{j=1}^{a_i-1}\frac{a_i-j}
{a_i}t_{i,j}\frac{\p}{\p t_{i,j}}
+\chi_A\frac{\p}{\p t_{\mu_A}}.
\end{equation}
\end{thm}
\begin{pf}
We use the following relation between the product $\circ$ 
and the intersection form $I_M$:
\begin{lem}\label{lem:for reconstruction}
Denote by $\Gamma^{ij}_{k}$ 
the contravariant components of the Levi--Civita connection for 
the intersection form $I_M$. 
Then one has
\begin{equation}
\Gamma^{ij}_{k}=d_{j} \ C^{\alpha \beta}_{\gamma},
\end{equation}
where $d$ is the dimension of the Frobenius manifold   
and $d_{j}$ is a rational number defined by $E(t_{j})=d_{j}t_{j}$,
$C^{ij}_{k}:=\sum_{a,b=1}^{\mu_A}\eta^{ai}\eta^{bj}C_{abk}$ and
$C_{ijk}:=\p_{i}\p_{j}\p_{k}\F_M$.  
\end{lem}
\begin{pf}
See Lemma 3.4 in \cite{d:1} and apply $d=1$. 
\qed
\end{pf}
One sees $C^{ij}_{k}$ can be reconstructed from
the intersection form $I_M$ if $d_{j}\ne 0$.
Since $d=1, d_{\beta}\neq 0$ if and only if $\beta=\mu_A$. 
However, we have 
\begin{equation*}
C^{i\mu_A}_{k}=\sum_{a,b=1}^{\mu_A}\eta^{ai}\eta^{b\mu_A}C_{abk}
=\delta^{i}_{k}, 
\end{equation*}
where $\delta^{i}_{k}$ is Kronecker's delta.
Therefore, $C_{ijk}$ and hence the Frobenius potential $\F_M$ 
can be reconstructed from the intersection form $I_M$
by Lemma~\ref{lem:for reconstruction}.
\qed
\end{pf}

\subsection{Isomorphism of Frobenius manifolds}

\begin{cor}\label{cor:sing and weyl}
Assume that $\chi_A>0$.
There exists an isomorphism of Frobenius manifolds between 
the one constructed from the invariant thory of extended affine Weyl group 
$\widehat{W}_A$ and the one constructed from the pair $(f_A,\zeta_A)$. 
\end{cor}
\begin{pf}
Corollary \ref{cor:sing and weyl} immediately follows from Theorem \ref{second}, 
Theorem \ref{calc: intersection form 2} and Theorem \ref{dz:thm}.
\qed
\end{pf}


\end{document}